\documentclass{amsart}
\usepackage{amsfonts, amssymb, verbatim}

\newtheorem{theorem}{Theorem}
\newtheorem{lemma}[theorem]{Lemma}
\newtheorem{proposition}[theorem]{Proposition}
\newtheorem{corollary}[theorem]{Corollary}

\theoremstyle{definition}
\newtheorem{definition}[theorem]{Definition}
\newtheorem{example}[theorem]{Example}

\theoremstyle{remark}
\newtheorem{remark}[theorem]{Remark}

\newtheorem{question}[theorem]{Question}

\numberwithin{equation}{section}
\numberwithin{theorem}{section}

\newcommand{\thref}[1]{Theorem~\ref{#1}}
\newcommand{\prref}[1]{Proposition~\ref{#1}}
\newcommand{\leref}[1]{Lemma~\ref{#1}}
\newcommand{\coref}[1]{Corollary~\ref{#1}}

\newcommand{\deref}[1]{Definition~\ref{#1}}
\newcommand{\exref}[1]{Example~\ref{#1}}

\newcommand{\reref}[1]{Remark~\ref{#1}}

\newcommand{\seref}[1]{Section~\ref{#1}}

\def\tt{\otimes}                               

\def\wti{\widetilde}

\def\ov{\overline}

\def\d{\partial}
\def\smash{\,\sharp\,}

\def\tint{{\textstyle\int}}

\def\otto{\leftrightarrow}                     
\def\isoto{\xrightarrow{\sim}}                 
\def\st{\; | \;}                               
\newcommand\nop[1]{{\rm{:}}\,{#1}\,{\rm{:}}}   
\def\vac{|0\rangle}                            

\def\op{{\mathrm{op}}}

\def\di{{\mathrm{d}}}
\def\CC{\mathbb{C}}       
\def\ZZ{\mathbb{Z}}       
\def\al{\alpha}                         

\def\Ga{\Gamma}
\def\de{\delta}
\def\De{\Delta}

\def\la{\lambda}

\def\A{\mathcal{A}}

\def\Mat{{\mathrm{Mat}}}
\def\glf{{\mathit{g{\ell}f}}}
\def\Zh{{\mathrm{Zh}}}         
\def\tens{\mathcal{T}}     
\def\univ{\mathcal{U}}     

\DeclareMathOperator{\Span}{span}

\DeclareMathOperator{\Res}{Res}
\DeclareMathOperator{\Ind}{Ind}

\DeclareMathOperator{\End}{End}

\DeclareMathOperator{\Lie}{Lie}

\newcommand{\alphaparenlist}{%
  \renewcommand{\theenumi}{\alph{enumi}}%
  \renewcommand{\labelenumi}{(\theenumi)}%
}

\newcommand{\romanparenlist}{%
  \renewcommand{\theenumi}{\roman{enumi}}%
  \renewcommand{\labelenumi}{(\theenumi)}%
}

\begin{document}

\title[Field Algebras]
{Field Algebras}

\author{Bojko Bakalov}
\address{Department of Mathematics, University of California,
Berkeley, CA 94720, USA}
\email{bakalov@math.berkeley.edu}
\thanks{The first author was supported by the Miller Institute
for Basic Research in Science.}

\author{Victor G.~Kac}
\address{Department of Mathematics,
MIT, Cambridge, MA 02139, USA}
\email{kac@math.mit.edu}
\thanks{The second author was supported in part by NSF grant DMS-9970007.}

\date{April 23, 2002. Revised May 6, 2002}

\begin{abstract}
A field algebra is a ``non-commutative'' generalization of a vertex algebra.  
In this paper we develop foundations of the theory of field algebras.
\end{abstract}

\dedicatory{Dedicated to Ernest Borisovich Vinberg on the occasion
of his 65th birthday.}

\maketitle


\setcounter{section}{-1}
\section{Introduction}
\label{sec:intro}

Roughly speaking, the notion of a vertex algebra \cite{B1} 
is a generalization of the notion of a unital commutative associative algebra
where the multiplication depends on a parameter.
(In fact, in \cite{B3} vertex algebras are described as
``singular'' commutative associative rings in a certain category.)
More precisely, an operator of left multiplication on a vertex algebra $V$
by an element $a \in V$ is a \emph{field} 
$Y(a,z) = \sum_{n\in\ZZ} a_{(n)} z^{-n-1}$,
where $a_{(n)} \in \End V$, and one requires that
$Y(a,z)b$ is a Laurent series in $z$ for any two
elements $a,b \in V$.
The role of a unit element of an algebra is played in this context
by a \emph{vacuum vector} $\vac \in V$, satisfying
\begin{list}{}{}
\item (\emph{vacuum axioms}) \quad $Y(\vac ,z) = I_V$, 
$Y(a,z) \vac = e^{zT} a$,
\end{list}
where $I_V\in\End V$ is the identity operator and $T\in\End V$.
A linear map $a \mapsto Y(a,z)$ of a vector space $V$ 
with a vacuum vector $\vac$ to the
space of $\End V$-valued fields satisfying 
\begin{list}{}{}
\item (\emph{translation invariance}) 
\quad $[T,Y(a,z)] = Y(Ta, z) = \partial_z Y(a,z)$,
\end{list}
is called a \emph{state--field correspondence}.
This notion is an analogue of a unital algebra.

For example, if $V$ is an ordinary algebra with a unit element $\vac$
and $T$ is a derivation of $V$, then the formula
  \begin{equation}\label{eq:0.1}
Y(a,z) b = (e^{zT}a) b , \qquad a,b \in V \,
  \end{equation}
defines a state--field correspondence
(and it is easy to show that all of them with the property that
$Y(a,z)$ is a formal power series in $z$, are thus obtained).

Furthermore, the associativity property of the algebra $V$
is equivalent to the following property of the state--field correspondence
\eqref{eq:0.1}:
  \begin{equation}\label{eq:0.2}
Y(Y(a,z)b,-w)c = Y(a,z-w)Y(b,-w)c \, , \qquad a,b,c \in V \, ,
  \end{equation}
and the commutativity property of $V$ is equivalent to:
  \begin{equation}\label{eq:0.3}
Y(a,z)b = e^{zT} \bigl( Y(b,-z)a \bigr) , \qquad a,b \in V \, .
  \end{equation}

It turns out that for a general vertex algebra $V$ identity \eqref{eq:0.2}
holds only after one multiplies both sides by $(z-w)^N$ where $N$
is sufficiently large (depending on $a,b,c$):
  \begin{equation}\label{eq:0.4}
(z-w)^N Y(Y(a,z)b,-w)c = (z-w)^N
     Y(a,z-w)Y(b,-w)c, \qquad N \gg 0 \, .
  \end{equation}
One of the equivalent definitions of a \emph{vertex algebra} is that it is a
state--field correspondence satisfying \eqref{eq:0.3} and \eqref{eq:0.4}
(see \thref{thm:vert}).

A \emph{field algebra} is a state--field correspondence satisfying 
only the associativity property \eqref{eq:0.4}.
We believe that this is the right analogue of a unital associative algebra.
In the present paper we are making the first steps towards
a general theory of field algebras.

The trivial examples of field algebras are provided by 
state--field correspondences \eqref{eq:0.1}: this is a
field algebra if and only if the underlying algebra $V$ is associative.
The simplest examples of non-trivial field algebras are tensor products
of field algebras \eqref{eq:0.1} with vertex algebras.
A special case of this is the algebra of matrices with entries
in a vertex algebra.
Other examples are provided by a smash product of a vertex algebra 
and a group of its automorphisms.
Thus, many important ring-theoretic constructions
with vertex algebras become possible in the framework of field algebras.

One of our main results is the construction of a canonical structure
of a field algebra in a tensor algebra $\tens(R)$
over a Lie (even Leibniz) conformal algebra $R$ (\thref{thm:t.1}).
Imposing relation \eqref{eq:0.3}
on $\tens(R)$ gives the enveloping vertex algebra $\univ(R)$ of $R$
(cf.\ \cite{K, GMS}).

We also establish the field algebra analogues of the density and
duality theorems in the representation theory of associative algebras
(see Theorems~\ref{th:5.1} and \ref{th:5.2})
and discuss the Zhu algebra construction \cite{Z} in the framework
of field algebras.

Note that the ``field algebras'' considered in \cite[Sec.~4.11]{K} 
are defined by a stronger than associativity axiom; we call them
strong field algebras in the present paper.
Surprisingly, it turns out that they are almost the same
as vertex algebras (see \thref{thm:sf.1}),
although the ``trivial'' field algebras \eqref{eq:0.1} are
automatically strong field algebras.
In the present paper the results of \cite[Sec.~4.11]{K}  on
field algebras are corrected.

The first examples of non-trivial field algebras (i.e., different from
\eqref{eq:0.1}) that are not vertex algebras were constructed in \cite{EK}. 
The ``quantum vertex algebras'' of \cite{EK} are 
field algebras satisfying in addition a certain ``braided commutativity''
generalizing \eqref{eq:0.3}. 
The relation of the present work to the paper \cite{EK}, and also to 
\cite{B3} and \cite{FR}, will be discussed in a subsequent paper.

\section{State--Field Correspondence}
\label{sec:sf}

Let $V$ be a vector space (referred to as the 
\emph{space of states}).  Recall that a (End~$V$-valued) \emph{field} 
is an expression of the form
\begin{equation*}
  a(z) = \sum_{n \in \ZZ} a_{(n)} z^{-n-1} \, , 
\end{equation*}
where $z$ is an indeterminate, $a_{(n)} \in \End V$, and for each 
$v \in V$ one has:
\begin{equation*}
  a_{(n)} v =0 \hbox{  for  } n \gg 0 \, , 
\end{equation*}
i.e., $a(z) v$ is a Laurent series in $z$.

Denote by $\glf(V)$ the space of all $\End V$-valued fields.  
For each $n \in \ZZ$ one defines the \emph{$n$-th product} 
of fields $a(z)$ and $b(z)$ by the following formula:
\begin{equation}
  \label{eq:1.1}
  a(z)_{(n)} b(z) = \Res_x \bigl( a(x)b(z) i_{x,z} (x-z)^n -
  b(z) a(x) i_{z,x} (x-z)^n \bigr) \, .
\end{equation}
Here $i_{x,z}$ (respectively $i_{z,x}$) stands for the expansion in the 
domain $|x| > |z|$ (respectively $|z| > |x|$):
\begin{equation*}
i_{x,z} (x-z)^n = \sum_{j=0}^\infty \binom{n}{j} x^{n-j} (-z)^j,
\end{equation*}
while
\begin{equation*}
i_{z,x} (x-z)^n = \sum_{j=0}^\infty \binom{n}{j} x^j (-z)^{n-j}.
\end{equation*}
%
%
It is easy to see that the space of fields $\glf(V)$ is closed under 
all $n$-th products and also under the derivation by the indeterminate.

Recall that
\begin{equation}
  \label{eq:1.2}
\de(x-z) = (i_{x,z} - i_{z,x}) (x-z)^{-1} 
= \sum_{j\in\ZZ} x^j z^{-j-1} 
\end{equation}
is the formal delta-function, characterized by the property:
\begin{equation*}
\Res_x a(x) \de(x-z) = a(z) \qquad \text{for }\; a(x) \in \glf(V).
\end{equation*}
Formula (\ref{eq:1.1}) is equivalent to the following two 
formulas for $n \in \ZZ_+$:
\begin{align*}
  a(z)_{(n)} b(z) &= \Res_x \, [a(x),b(z)] (x-z)^n \, , \\
  a(z)_{(-n-1)} b(z) &= \nop{ \partial^n_z a(z) \, b(z) } /n! \, .
\end{align*}
Here $:\;\;:$ stands for the \emph{normal ordered product} of fields defined by
\begin{equation}\label{eq:nop}
 \nop{ a(z) b(z) } = a(z)_+ b(z) + b(z)a(z)_- \, ,
\end{equation}
where
\begin{equation*}
  a(z)_+ = \sum_{j \leq -1} a_{(j)} z^{-j-1} \, , \quad
  a (z)_- = \sum_{j \geq 0} a_{(j)} z^{-j-1} \, .
\end{equation*}

\begin{definition}
  \label{def:1.1}
Let $(V, \vac)$ be a \emph{pointed vector space}, i.e., a
vector space $V$ with a fixed non-zero vector $\vac \in V$, 
which will be referred to as the \emph{vacuum vector}.
A \emph{state--field correspondence} 
is a linear map 
\begin{equation*}
Y\colon V \to \glf(V), \quad
a \mapsto Y(a,z) = \sum_{n\in\ZZ} a_{(n)} z^{-n-1} ,
\end{equation*}
such that following axioms hold:
\begin{list}{}{}
\item (\emph{vacuum axioms}) \quad $Y(\vac ,z) = I_V$, 
$Y(a,z) \vac = a + T(a)z + \cdots \in V[[z]]$,
\end{list}
where $I_V\in\End V$ is the identity operator, $T\in\End V$,
%
\begin{list}{}{}
\item (\emph{translation invariance}) 
\quad $[T,Y(a,z)] = Y(Ta, z) = \partial_z Y(a,z)$.
\end{list}
The linear operator $T$ on $V$ is called the \emph{translation operator}.  

\end{definition}

\begin{example}\label{ex:1.1}
  Let $V$ be a unital algebra with a unit element $\vac$, and let 
  $T$ be a derivation of $V$.  Then
\begin{equation*}
    Y(a,z) b = (e^{zT}a) b , \qquad a,b \in V \, ,
\end{equation*}
 is a state--field correspondence.  In fact, all state--field correspondences
for which the fields $Y(a,z)$ are formal power series in $z$ 
are obtained in this way.
We will call such a state--field correspondence \emph{trivial}.
\end{example}

\begin{example}\label{ex:1.3}
Given two state--field correspondences $(V_i, \vac_i, Y_i)$, $i=1,2$,
one defines their tensor product 
$(V = V_1 \tt V_2, \vac = \vac_1 \tt \vac_2, Y)$,
where $Y(a_1 \tt a_2, z) = Y(a_1,z) \tt Y(a_2,z)$.
This is again a state--field correspondence, the 
translation operator being $T = T_1 \tt I + I \tt T_2$.
\end{example}

\begin{example}\label{ex:1.4}
A special case of \exref{ex:1.3} is the tensor product of a
state--field correspondence $(V, \vac, Y)$ and a unital algebra
$(A,1)$, which is viewed as a trivial state--field correspondence
with $T=0$. In particular, the space $\Mat_N(V)$ of $N$ by $N$ matrices
with entries in $V$ has a structure of a state--field correspondence.
\end{example}

\begin{example}\label{ex:1.5}
A generalization of \exref{ex:1.4} is the \emph{smash product}
$V \smash \Ga$, where $\Ga$ is a group of automorphisms of the
state--field correspondence $(V, \vac, Y)$. We define the
space of states to be $V \tt \CC[\Ga]$, the vacuum vector to be
$\vac \tt 1$, and let
\begin{equation*}
    Y(a \tt g, z) (b \tt h) = Y(a,z) (g b) \tt gh \, ,
\qquad a,b \in V , \; g,h \in \Ga \, .
\end{equation*}
When $\Ga$ is finite, one has the following very useful formula 
(cf.\ \cite{MS}):
\begin{equation}\label{eq:1.3}
V^\Ga \simeq \{ e_{(-1)} (v_{(-1)} e) \st v \in V \smash \Ga \} 
\subset V \smash \Ga \, ,
\qquad \text{where }\; e = \frac1{|\Ga|} \sum_{g\in\Ga} \vac \tt g \, .
\end{equation}
Here $V^\Ga = \{ a \in V \st g a = a \text{ for all } g \in \Ga \}$
denotes the subalgebra of invariants of $\Ga$ in~$V$. 
The isomorphism \eqref{eq:1.3} is given by 
$a \mapsto \frac1{|\Ga|} \sum_{g\in\Ga} a \tt g$
for $a \in V^\Ga$.
\end{example}

One defines in the obvious way the notions of subalgebras, ideals
and homomorphisms of state--field correspondences. For example,
a \emph{subalgebra} of a state--field correspondence $(V, \vac, Y)$ 
is a subspace $U$ of $V$ containing $\vac$ and such that
$a_{(n)} b \in U$ for all $n\in\ZZ$ if $a,b\in U$.
A \emph{left ideal} is a $T$-invariant subspace $I$ of $V$ such that
$a_{(n)} b \in I$ for all $n\in\ZZ$ if $a\in V$, $b\in I$.
If $I$ is a two-sided ideal of $V$, there is a natural state--field 
correspondence with a space of states $V/I$.

\begin{proposition}\label{prop:1.sf}
In terms of Fourier coefficients $a_{(n)}$ of\/ $Y$,
the definition of a state--field correspondence
can be reformulated as follows{\rm:}
\begin{list}{}{}
\item \emph{(local finiteness)} \;
$a_{(N)} b = 0$ \; for $a,b\in V$, $N\gg0$,

\item \emph{(weak vacuum axiom)} \;
$a_{(-1)} \vac = \vac_{(-1)} a = a$, 
\item \emph{(translation invariance)} \,
$[T, a_{(n)}] = (Ta)_{(n)} = -n a_{(n-1)}$ \; for $n\in \ZZ$, $a\in V$.
\end{list}
\end{proposition}
\begin{proof}
We have to check that 
$\vac_{(n)} a = 0$ for $n\ne -1$, $Ta = a_{(-2)} \vac$, 
and $a_{(n)} \vac = 0$ for $n\ge 0$, $a\in V$.
First, $Ta = (Ta)_{(-1)} \vac = a_{(-2)} \vac$.
Next, we have $T\vac = T(\vac_{(-1)}\vac) = 2T\vac$, hence $T\vac=0$.
Then for every $n\ne0$, $(T\vac)_{(n)} = -n \vac_{(n-1)}$
shows that $\vac_{(n-1)} = 0$.
Similarly, $T(a_{(n)} \vac) = -n a_{(n-1)} \vac$ implies 
$a_{(n)} \vac = 0$ for $n\ge 0$, because $a_{(N)} \vac = 0$ for $N\gg0$.
\end{proof}

\begin{proposition}
  \label{prop:1.1}
  
Any state--field correspondence $Y$ has the following properties{\rm:}

\alphaparenlist
\begin{enumerate}
\item 
$Y(a,z) \vac = e^{zT} a$,

\item 
$e^{wT} Y(a,z) e^{-wT} = Y(e^{wT} a, z) = i_{z,w} Y(a, z+w)$,

\item 
$(Y(a,z)_{(n)}Y(b,z)) \vac = Y(a_{(n)}b, z) \vac $.


\end{enumerate}
\end{proposition}
\begin{proof}
  
See \cite[Proposition~4.1]{K}. 
\end{proof}

\begin{proposition}
  \label{prop:1.2}
Given a state--field correspondence $Y$, define 
\begin{equation}\label{eq:1.op} 
Y^\op(a,z)b = e^{zT}Y(b,-z)a.
\end{equation}
Then $Y^\op$ is also a state--field correspondence,
and $(Y^\op)^\op = Y$.
\end{proposition}
\begin{proof}
Straightforward.
\end{proof}

\begin{definition}
  \label{def:1.3}
%
%
%
The state--field correspondence $Y^\op$ defined by \eqref{eq:1.op} is 
called the \emph{opposite to} $Y$.
\end{definition}

At the end of this section, we study the notion of grading for
state--field correspondences.

\begin{definition}
  \label{def:g.1}
A state--field correspondence $Y$ on a pointed vector space
$(V,\vac)$ is called \emph{graded}
if there is a diagonalizable operator $H \in \End V$ satisfying
\begin{equation}\label{eq:g.1}
[H, Y(a,z)] = Y(Ha,z) + z\d_z Y(a,z) \, ,
\qquad a \in V \, .
\end{equation}
It is called $\ZZ_+$-\emph{graded} if 
all eigenvalues of $H$ are non-negative integers.
The operator $H$ is called a \emph{Hamiltonian} of $V$.
\end{definition}

For a homogeneous element $a \in V$, we denote its degree by
$\De_a$: $Ha = \De_a a$.
In terms of modes, equation \eqref{eq:g.1} is equivalent to:
\begin{equation}\label{eq:g.2}
[H, a_{(n)}] = (\De_a-1-n) a_{(n)}  \, ,
\qquad a \in V, \;\; Ha = \De_a a, \;\; n \in \ZZ \, .
\end{equation}
This implies that for homogeneous $a,b \in V$ one has:
\begin{equation}\label{eq:g.2a}
\De_{a_{(n)} b} = \De_a+\De_b-n-1 \, , \qquad n \in \ZZ \, .
\end{equation}
One can easily show that $H\vac = 0$ and $[H,T]=T$.
The latter is equivalent to: $\De_{Ta} = \De_a+1$ for homogeneous $a \in V$
(and is a special case of \eqref{eq:g.2a} since $Ta=a_{-2}\vac$). 

\begin{proposition}\label{prop:g.1}
If\/ $Y$ is a graded state--field correspondence,
then its opposite $Y^\op$ is also graded,
with the same Hamiltonian $H$.
\end{proposition}
\begin{proof}
Easy exercise, using \eqref{eq:1.op}, \eqref{eq:g.1} and 
$H e^{zT} = e^{zT} (H+zT)$.
\end{proof}

\begin{remark}\label{rem:g.1}
A field $\phi(z) \in \glf(V)$ is said to have \emph{conformal dimension}
$\De$ if
\begin{equation*}\label{eq:g.6}
[H, \phi(z)] = (\De + z\d_z) \phi(z) \, .
\end{equation*}
Then for each homogeneous $a \in V$, $Y(a,z)$ is a field of 
conformal dimension $\De_a$ (see \eqref{eq:g.1}). 

Let $\phi(z), \phi'(z) \in \glf(V)$ be two fields of 
conformal dimensions $\De$ and $\De'$, respectively.
Then $\d_z \phi(z)$ has conformal dimension $\De+1$, and
for each $n\in\ZZ$, the $n$-th product $\phi(z)_{(n)} \phi'(z)$
has conformal dimension $\De+\De'-n-1$.

\end{remark}

If $Y$ is a $\ZZ_+$-graded state--field correspondence,
we may shift subscripts as follows:
\begin{equation}\label{eq:g.3}
a_n = a_{(n+\De_a-1)} \, , \quad
a_{(n)} = a_{n-\De_a+1} \qquad\text{if} \;\; Ha = \De_a a \, ,
\end{equation}
so that
\begin{equation}\label{eq:g.4}
[H, a_n] = -n \, a_n  \, ,
\qquad a \in V, \; n \in \ZZ \, .
\end{equation}
It is easy to see that:
\begin{equation}\label{eq:g.5}
((T+H)a)_0 = 0 \qquad\text{for all} \;\; a \in V \, .
\end{equation}

\section{Locality}
\label{sec:loc}

  A pair of fields $a(z),b(z) \in \glf(V)$ is called \emph{local} if 
\begin{equation*}
(z-w)^N [a(z),b(w)]=0 \quad\text{for } N\gg0.
\end{equation*}
It is called \emph{local on $v \in V$} if
\begin{equation*}
(z-w)^N [a(z),b(w)] \, v=0 \quad\text{for } N\gg0.
\end{equation*}

\begin{proposition}
\label{prop:loc}
A pair of fields $(a(z),b(z))$ is local if and only if
\begin{equation*}
  \label{eq:loc1}
  [a(z),b(w)] = \sum_{\substack{j \geq 0\\ \text{\rm{finite}}}} 
  \bigl( a(w)_{(j)} b(w) \bigr) \,
     \partial^j_w \delta (z-w)/j! \, .
\end{equation*}
It is local on $v$ iff
\begin{equation*}
  \label{eq:loc2}
  [a(z),b(w)] \, v = \sum_{\substack{j \geq 0\\ \text{\rm{finite}}}} 
  \bigl( a(w)_{(j)} b(w) \bigr) v \,
     \partial^j_w \delta (z-w)/j! \, .
\end{equation*}
\end{proposition}
\begin{proof}
Follows from \cite[Corollary~2.2]{K}.  
\end{proof}

A pair of fields $(a(z),b(z))$ is called \emph{weakly local} if
\begin{equation*}
\Res_z (z-w)^N [a(z),b(w)]=0 \quad\text{for } N\gg0.
\end{equation*}
Note that the weak locality of the pair $(a(z),b(z))$ means that 
$a(z)_{(N)} b(z) =0$ for $N \gg 0$.  

\begin{remark}
\label{rem:wl}
If a pair of fields $(a(z),b(z))$ is local, then the pair 
$(b(z),a(z))$ is also local. 
However, the weak locality of a pair
$(a,b)$ does not imply the weak locality of the pair $(b,a)$.
Indeed, let $a(z)$ be the free boson, so that 
$[a_{(m)},a_{(n)}] = m \, \de_{m,-n}$, and take 
$b(z) = a_{(1)} z^{-1}$.
Then $[a(z), b(w)] = -w^{-1}$; therefore, 
$a(z)_{(n)} b(z) =0$, $b(z)_{(n)} a(z) = (-z)^n$ for all $n\ge0$.
%
\end{remark}

  A collection of fields is called \emph{local} (respectively 
\emph{local on $v \in V$}, respectively \emph{weakly local}) 
if each pair of fields from this collection is 
local (respectively local on $v \in V$, respectively weakly local).

A state--field correspondence $Y$ is called \emph{local} 
(respectively \emph{weakly local}) 
if the collection of fields $\{ Y(a,z) \}_{a \in V}$ is local 
(respectively weakly local).

\begin{example}
  \label{ex:1.2}
Let $U$ be a vector space and let $V \subset \glf(U)$ be 
a weakly local space of
$\End U$-valued fields in the variable $x$.  
Let $\vac =I_U \in V$.
Then the following formula defines a state--field correspondence:
\begin{equation*}
  Y(a(x), z) b(x) = \sum_{n \in \ZZ} a(x)_{(n)} b(x) \, z^{-n-1} \, ,
  \qquad a(x),b(x) \in V,
\end{equation*}
the translation operator being $T= \partial_x$.  
\end{example}

The following proposition generalizes the Uniqueness theorem
from \cite[Sec.~4.4]{K}.

\begin{proposition}
\label{prop:2.1}
Let $Y$ be a state--field correspondence on
a pointed vector space $(V, \vac)$.  
Let $B_i (z) \in \glf(V)$, $i=1,2$, be such that
\romanparenlist
\begin{enumerate}
\item 
$B_1 (z) \vac = B_2 (z) \vac$,

\item 
all pairs $(Y(a,z),B_i (z))$ are local on $\vac$.

\end{enumerate}
Then $B_1 (z) = B_2 (z)$.
\end{proposition}
\begin{proof}
Let $B(z) = B_1 (z) - B_2 (z)$.  Then $B(z) \vac =0$ and all pairs 
$(Y(a,z),B(z))$ are local on $\vac$. We have $(z-w)^N [B(z),Y(a,w)] \vac =0$, 
so $(z-w)^N B(z) Y(a,w) \vac$ $=0$ for $N\gg0$.
Using Proposition~\ref{prop:1.1}a, we get
$(z-w)^N B(z) e^{wT} a =0$.
Letting $w=0$, we obtain $B(z) a =0$, $a \in V$.
\end{proof}

The following is a generalization of Dong's lemma 
(see \cite[Lemma~3.2]{K}).

\begin{lemma}\label{ldong}
Let $a,b,c \in\glf(V)$ be three fields.

\alphaparenlist
\begin{enumerate}
\item 
If the pair $(a,b)$ is weakly local and the pairs
$(a,c)$, $(b,c)$ are local, then 
$(a_{(n)} b,c)$ is local for any $n\in\ZZ$.

\item 
If the pairs $(a,b)$ and $(a,c)$ are weakly local and the pair
$(b,c)$ is local, then 
$(a_{(n)} b,c)$ is local for any $n\ge0$.

\item 
If the pairs $(a,b)$, $(a,c)$, $(b,c)$ are weakly local, then
$(a_{(n)} b,c)$ is weakly local for any $n\ge0$.

\item 
If the pairs $(a,b)$, $(a,c)$, $(b,c)$ are weakly local, then
$(a, b_{(n)} c)$ is weakly local for any $n\in\ZZ$.

\end{enumerate}

\end{lemma}
\begin{proof}
Same as in \cite[Sec.~3.2]{K}.
\end{proof}

\begin{remark}
\label{rem:w2}
In general, it is not true that if all pairs $(a,b)$, $(a,c)$, $(b,c)$ 
are weakly local, then $(a_{(-1)} b, c)$ is weakly local as well.
For example, let $a(z)$ be the free boson (see \reref{rem:wl}),
and take $b(z)=a(z)$, $c(z)=a(z)_+$.
Then $[a(z)_-, a(w)_+] = i_{z,w} (z-w)^{-2}$ and $[a(z)_+, a(w)_+] = 0$;
hence $(a,a)$ and $(a,a_+)$ are weakly local.
We have: $[\nop{ a(z)^2 }, a(w)_+] = [a(z)_+ a(z) + a(z) a(z)_-, a(w)_+] =
2 a(z) \, i_{z,w} (z-w)^{-2}$, which shows that the pair
$(\nop{ a^2 }, a_+)$ is not weakly local.
\end{remark}

\begin{lemma}\label{lem:2.2}
Let $X$ and $Y$ be two state--field correspondences on
a pointed vector space $(V, \vac)$,
and let $a,b,c \in V$ be such that

\romanparenlist
\begin{enumerate}
\item 
  $(Y(a,z),Y(b,z))$ is a weakly local pair,

\item 
  the pair $(Y(a,z), X(c,z))$ {\rm(}respectively, $(Y(b,z),X(c,z))${\rm)}
  is local on $\vac$ and $b$ {\rm(}respectively, on $\vac$ and $a${\rm)}.
\end{enumerate}
 Then the pair $(Y(a,z)_{(n)} Y(b,z), X(c,z))$ 
is local on the vacuum vector for all $n \in \ZZ$.

\end{lemma}
\begin{proof}
Let us write $a(z) = Y(a,z)$, $b(z) = Y(b,z)$, $c(z) = X(c,z)$ for short.
By the (weak) locality, there exists a number $N\ge0$ 
such that for all $r\ge N$:
\begin{align*}
\Res_{z_1} (z_1-z_2)^r &[a(z_1),b(z_2)] = 0,
\\
(z_1-z_3)^r &[a(z_1),c(z_3)] \, v = 0, \qquad v=\vac \text{ or } b \, ,
\\
(z_2-z_3)^r &[b(z_2),c(z_3)] \, v = 0, \qquad v=\vac \text{ or } a \, .
\end{align*}
Let us apply the differential operator
$e^{-z_2 ( \d_{z_1} + \d_{z_3} )}$
to the left-hand side of the second equation 
for $v=b$.
By Taylor's formula, the result is:
\begin{equation*}
(z_1-z_3)^r i_{z_1,z_2} i_{z_3,z_2} 
[a(z_1-z_2),c(z_3-z_2)] \, b \, .
\end{equation*}
Using \prref{prop:1.1}, this is equal to:
\begin{align*}
& (z_1-z_3)^r [e^{-z_2 T} a(z_1) e^{z_2 T}, e^{-z_2 T} c(z_3) e^{z_2 T}] \, b
\\
= & (z_1-z_3)^r e^{-z_2 T} [a(z_1),c(z_3)] \, e^{z_2 T} b 
\\
= & (z_1-z_3)^r e^{-z_2 T} [a(z_1),c(z_3)] \, b(z_2)\vac \, .
\end{align*}
Applying $e^{z_2 T}$ to this, we obtain:
\begin{equation*}
(z_1-z_3)^r [a(z_1),c(z_3)] \, b(z_2)\vac = 0, \qquad r\ge N.
\end{equation*}
Similarly, we have:
\begin{equation*}
(z_2-z_3)^r [b(z_2),c(z_3)] \, a(z_1)\vac = 0, \qquad r\ge N.
\end{equation*}
The rest of the proof is as in \cite[Lemma~3.2]{K}.
\end{proof}

In the proof of \leref{lem:2.2}, we proved the following result
which will be needed later.

\begin{lemma}\label{lem:2.2a}
Let $X$ and $Y$ be two state--field correspondences,
and let $a,b,c \in V$ be such that
\begin{equation*}
(z-w)^N [Y(a,z), X(c,z)] \, b = 0 
\end{equation*}
for some $N\ge0$. Then
\begin{equation*}
(z-w)^N [Y(a,z), X(c,z)] \, T b = 0 \, .
\end{equation*}
\end{lemma}

\begin{proposition}
\label{prop:2.3}
Let $X$ and $Y$ be two state--field correspondences for $(V, \vac)$,
 such that

\romanparenlist
\begin{enumerate}
\item 
  $Y$ is weakly local,

\item 
  all  pairs $(Y(a,z), X(b,z))$  are local on any $v \in V$.
\end{enumerate}
Then $Y(a,z)_{(n)} Y(b,z)= Y(a_{(n)} b,z)$ for all $n \in \ZZ$, $a,b \in V$.

\end{proposition}
\begin{proof}
Let $B_1 (z) = Y(a,z)_{(n)} Y(b,z)$, $B_2 (z) = Y(a_{(n)}b,z)$.  
Due to Proposition~\ref{prop:1.1}c, the fields $B_i (z)$ satisfy condition~(i) 
of Proposition~\ref{prop:2.1}.  Due to \leref{lem:2.2}, 
condition~(ii) of  Proposition~\ref{prop:2.1} is satisfied as well.  
Hence $B_1 (z) = B_2 (z)$.
\end{proof}

\begin{proposition}\label{prop:2.8}
Let $X$ and $Y$ be two state--field correspondences. Then
all pairs $(Y(a,z), X(b,z))$ are local on the vacuum vector
if and only if $X=Y^\op$.
\end{proposition}
\begin{proof}
Assume that $(z-w)^N Y(a,z) X(b,w) \vac =
(z-w)^N X(b,w) Y(a,z) \vac$ for $N\gg0$.
Using \prref{prop:1.1}, this is equivalent to:
$(z-w)^N e^{wT} i_{z,w} Y(a,z-w) b = (z-w)^N e^{zT} i_{w,z} X(b,w-z) a$.
For a sufficiently large $N$, both sides contain only non-negative
powers of $z-w$, and hence only non-negative powers of $z$ and $w$.
Putting $z=0$, we obtain 
$X(b,w) a = e^{wT} Y(a,-w) b$.

Conversely, if $X=Y^\op$, the above calculations show that
$[Y(a,z), X(b,w)] \vac =$~$0$ for all $a,b \in V$.
\end{proof}

\section{Field Algebras}
\label{sec:field}

Let $(V, \vac)$ be a pointed vector space and let $Y$ be a
state--field correspondence.
We say that $Y$ satisfies the \emph{$n$-th product axiom} if
for all $a,b\in V$ and $n\in\ZZ$
\begin{equation}
   \label{eq:3.4}    
Y(a_{(n)} b, z) = Y(a,z)_{(n)} Y(b,z) \, .
\end{equation}
We say that $Y$ satisfies the \emph{associativity axiom}
if for all $a,b,c\in V$
\begin{equation}
  \label{eq:3.9}
  (z-w)^N Y(Y(a,z)b,-w)c = (z-w)^N i_{z,w}
     Y(a,z-w)Y(b,-w)c, \qquad N \gg 0 \, .
\end{equation}

\begin{proposition}
\label{prop:ax}
Let $(V, \vac)$ be a pointed vector space and let $Y$ be a
state--field correspondence. 
Then{\rm:}

\alphaparenlist
\begin{enumerate}
\item 
$Y$ satisfies the $n$-th product axiom \eqref{eq:3.4} iff for all $a,b\in V$
\begin{equation}
  \label{eq:lr}
  [Y(a,z),Y^\op(b,w)] = \sum_{\substack{j \geq 0\\ \text{\rm{finite}}}}
  Y^\op(a_{(j)} b,w) \,
     \partial^j_w \delta (z-w)/j! \, ,
\end{equation}
where $Y(a,z) =\sum_{j\in\ZZ} a_{(j)}z^{-j-1}$
and $Y^\op$ is the opposite to $Y$ {\rm(}see \eqref{eq:1.op}{}{\rm)}.

\item 
$Y$ satisfies the associativity axiom \eqref{eq:3.9} iff all 
pairs $(Y(a,z),Y^\op(b,z))$ are local on each $v \in V$.
In particular, the $n$-th product axiom implies
the associativity axiom.

\end{enumerate}
\end{proposition}
\begin{proof}
Assume that $Y$ satisfies the $n$-th product axiom.
Replace $z$ with $-w$ in
(\ref{eq:3.4}), multiply both sides by $z^{-n-1}$, and sum over
$n \in \ZZ$.  Then after some manipulation, we get:
%
\begin{equation}\label{eq:3.8}
\begin{split}
  Y(Y&(a,z)b,-w)c \\
    &= i_{z,w} Y(a,z-w)Y(b,-w)c
    -Y(b,-w)\sum_{j \geq 0} (a_{(j)}c) \, \partial^j_w
    \delta (z-w)/j! \, .
\end{split}
\end{equation}
After applying $e^{wT}$ to both sides, the left-hand side
becomes $Y^\op(c,w)Y(a,z)b$. The first term in the right-hand side
becomes (using \prref{prop:1.1}):
\begin{equation*}
i_{z,w} e^{wT} Y(a,z-w)Y(b,-w)c = Y(a,z) e^{wT} Y(b,-w)c = 
Y(a,z) Y^\op(c,w) b,
\end{equation*}
while the second term becomes:
\begin{equation*}
\sum_{j \geq 0} Y^\op(a_{(j)}c,w)b \, \partial^j_w
    \delta (z-w)/j! \, .
\end{equation*}
Therefore \eqref{eq:3.8} is equivalent to \eqref{eq:lr}
(with $b$ replaced with $c$).
Since \eqref{eq:3.8} is the generating function of all
$n$-th products, we get statement (a).

Notice that, by the above argument, the locality of the pair
$(Y(a,z),Y^\op(c,z))$ on $b$ is equivalent to the associativity property
\eqref{eq:3.9}. This completes the proof.
\end{proof}

\begin{corollary}
  \label{cor:3.2}
Let $Y$ be a state--field correspondence for $(V, \vac)$.
Then for three elements $a,b,c \in V$, the collection of identities 
\begin{equation}\label{eq:lln}
Y(a_{(n)} b, z) c = \bigl( Y(a,z)_{(n)} Y(b,z) \bigr) c,
\qquad n\in\ZZ
\end{equation}
implies the locality of the pair $(Y(a,z),Y^\op(c,z))$ on $b$,
which in turn is equivalent to the associativity property \eqref{eq:3.9}.
If \eqref{eq:lln} holds for fixed $a,c$
and all $b \in V$, then the pair $(Y(a,z),Y^\op(c,z))$ is local.
\end{corollary}

\begin{definition}
  \label{def:3.1}
Let $(V,\vac)$ be a pointed vector space.
A \emph{field algebra} $(V,\vac,Y)$ is a state--field
correspondence $Y$ for $(V,\vac)$ satisfying 
the associativity axiom \eqref{eq:3.9}.
A \emph{strong field algebra} $(V,\vac,Y)$ is a state--field
correspondence $Y$ satisfying the $n$-th product axiom \eqref{eq:3.4}.
\end{definition}

Note that, by \prref{prop:ax}b, any strong field algebra is
a field algebra.

\begin{example}\label{ex:3.3}
Let $V$ be an associative algebra with a unit element $\vac$ and
let $T$ be a derivation of $V$.  Let $Y(a,z)b = (e^{zT}a) b$.
Then $(V,\vac,Y)$ is a strong field algebra. 
All field algebras for which all $Y(a,z)$ are formal power series in
$z$ are obtained in this way. We call such field algebras \emph{trivial}.
\end{example}

\begin{example}\label{ex:3.5}
The tensor product of two field algebras (cf.\ \exref{ex:1.3})
is again a field algebra. However, the tensor product of two 
strong field algebras is not necessarily a strong field algebra.
For example, the tensor product of a non-trivial strong field algebra
and a trivial field algebra (see \exref{ex:3.3}), where the underlying
associative algebra is non-commutative,
is not a strong field algebra.
\end{example}

\begin{example}\label{ex:3.6}
The smash product of a field algebra $V$ and a group $\Ga$ of its 
automorphisms (see \exref{ex:1.5}) is a field algebra. It is not 
a strong field algebra if $V$ is non-trivial and the action of $\Ga$
on $V$ is non-trivial.
\end{example}

\begin{theorem}\label{th:3.1}
\alphaparenlist
\begin{enumerate}
\item 
A field algebra $(V,\vac,Y)$ is the same as a state--field
correspondence $Y$ for $(V,\vac)$ such that there
exists a state--field correspondence $X$ for $(V,\vac)$,
having the property that all pairs $(Y(a,z), X(b,z))$ are local on
each $v \in V$. In this case $X=Y^\op$.

\item 
A strong field algebra is the same as a field algebra $(V,\vac,Y)$
for which the state--field correspondence $Y$ is weakly local.

\end{enumerate}
\end{theorem}
\begin{proof}
(a) Let $(V,\vac,Y)$ be a field algebra.  Then, by \prref{prop:ax}b,
$Y^\op$ is local with $Y$ on each $v \in V$.

Conversely, let $X$ be a state--field correspondence which is 
local with $Y$ on every vector. In particular, it is local on the vacuum.
Then, \prref{prop:2.8} implies that $X=Y^\op$, and 
again by \prref{prop:ax}b, $(V,\vac,Y)$ is a field algebra.

(b) If $(V,\vac,Y)$ is a strong field algebra, then it is also a 
field algebra. The weak locality of $Y$ follows from \eqref{eq:3.4},
because $a_{(n)} b = 0$ for $n\gg0$.

Conversely, if $(V,\vac,Y)$ is a field algebra with a weakly local $Y$,
the $n$-th product axiom \eqref{eq:3.4} follows from part (a) and
\prref{prop:2.3}.  
\end{proof}

\begin{remark}
  \label{rem:3.1}
If $(V,\vac,Y)$ is a field algebra,
then $(V,\vac,Y^\op)$ is also a field algebra, 
called the \emph{opposite} field algebra.
\end{remark}

\begin{corollary}\label{cor:fie2}
Let $(V,\vac,Y)$ be a strong field algebra.
Then for all $a,b\in V$ the pair $(Y(a,z),Y^\op(b,z))$ is local and
\begin{equation*}\label{eq:lrn}
Y(a,z)_{(n)} Y^\op(b,z) = Y^\op(a_{(n)} b,z), \qquad n\ge0.
\end{equation*}
\end{corollary}
\begin{proof}
This follows from \prref{prop:ax}a.
\end{proof}

\begin{proposition}\label{prop:fie}
  A strong field algebra is a vector space $V$ with a given vector
  $\vac$ and a linear map $Y\colon V \to \glf(V)$, $a \mapsto
  Y(a,z) =\sum_{n\in\ZZ} a_{(n)}z^{-n-1}$, satisfying for all $n \in
  \ZZ$, $a,b \in V$ the $n$-th product axiom \eqref{eq:3.4} and
  the weak vacuum axioms
\begin{equation}\label{eq:wvax}
\vac_{(n)} a = \delta_{n,-1}a, \quad a_{(-1)} \vac =a \, .
\end{equation}
%
\end{proposition}
\begin{proof}
The only thing that has to be checked is that such a map $Y$
is a state--field correspondence. 
Let $Ta = a_{(-2)} \vac$.  Letting $b=\vac$ in
$n$-th product axioms for $n \geq 1$, 
we see that the vacuum axioms hold.  Notice that, taking coefficient of
$w^{-k-1}$, the $n$-th product axiom gives the following
relation:
\begin{equation}\label{eq:3.10}
  (a_{(n)}b)_{(k)}c = \sum^{\infty}_{j=0} (-1)^j \binom{n}{j}
  \bigl( a_{(n-j)} (b_{(k+j)} c) - (-1)^n b_{(n+k-j)} (a_{(j)} c) \bigr) \, .
\end{equation}
Letting $k=-2$ and $c=\vac$ in (\ref{eq:3.10}) gives 
$[T,a_{(n)}] =-n \, a_{(n-1)}$.
Letting $n=-2$ and $b=\vac$ in (\ref{eq:3.10}), we get 
$(Ta)_{(k)} =-k \, a_{(k-1)}$.
Therefore, $Y$ is a state--field correspondence. 
\end{proof}

\begin{definition}
  \label{def:g.2}
A (strong) field algebra $(V,\vac,Y)$ is called \emph{graded},
respectively $\ZZ_+$-\emph{graded}, if the 
state--field correspondence $Y$ is (see \deref{def:g.1}).
\end{definition}
\begin{remark}\label{rem:g.2}
Let $(V,\vac,Y)$ be a $\ZZ_+$-graded strong field algebra, 
and let $a,b,c \in V$ be homogeneous elements. Then
$Y(a,z)_{(n)} Y(b,z) = 0$ for $n \ge \De_a+\De_b$,
and the associativity relation \eqref{eq:3.9} holds for
$N \ge \De_a+\De_c$.
This follows from the proof of \prref{prop:ax} and 
\eqref{eq:g.2a}. 
\end{remark}

\begin{question}
Is it true that a weakly local subspace $V$ of $\glf(U)$,
containing $I_U$, $\partial_x$-invariant and closed under all $n$-th products,
is a strong field algebra? (See \exref{ex:1.2}.)
\end{question}

\section{Conformal Algebras and Field Algebras}
\label{sec:conf}

Let $(V, \vac)$ be a pointed vector space and let $Y$ be a
state--field correspondence.
For $a,b\in V$, we define
their \emph{$\la$-product} by the formula
\begin{equation*}\label{eq:c.1}
a_\la b = \Res_z e^{\la z} Y(a,z)b 
= \sum_{\substack{n \geq 0\\ \text{finite}}} \la^n a_{(n)} b / n! \, .
\end{equation*}
We also have the $(-1)$-st product on $V$, which we denote as
\begin{equation*}\label{eq:c.2}
a.b = \Res_z z^{-1} Y(a,z)b = a_{(-1)} b .
\end{equation*}
The vacuum axioms for $Y$ imply:
\begin{equation}\label{eq:c.3}
\vac.a = a = a.\vac ,
\end{equation}
while the translation invariance axiom shows that:
\begin{equation}\label{eq:c.6}
T (a.b) = (Ta).b + a.(Tb)
\end{equation}
and
\begin{equation}\label{eq:c.4}
T (a_\la b) = (Ta)_\la b + a_\la (Tb) , \quad
(Ta)_\la b = -\la \, a_\la b  
\end{equation}
for all $a,b\in V$.

Conversely, if we are given a linear operator $T$,
a $\la$-product and a $.$-product on $V$,
satisfying the above properties \eqref{eq:c.3}--\eqref{eq:c.4},
we can reconstruct the state--field correspondence $Y$ by the formulas:
\begin{equation}\label{eq:c.5}
Y(a,z)_+ b = (e^{zT} a).b , \qquad
Y(a,z)_- b = (a_{-\d_z} b) (z^{-1}),
\end{equation}
where $Y(a,z) = Y(a,z)_+ + Y(a,z)_-$.
Notice that equations \eqref{eq:c.3}--\eqref{eq:c.4}
imply $T\vac = 0$ and $\vac_\la a = 0 = a_\la \vac$ for $a\in V$
(cf.\ \prref{prop:1.sf}).



A $\CC[T]$-module $V$, equipped with a linear map 
$V \tt V \to \CC[\la] \tt V$, $a\tt b \mapsto a_\la b$,
satisfying \eqref{eq:c.4} is called a ($\CC[T]$-)\emph{conformal algebra}
(cf.\ \cite[Sec.~2.7]{K}).
On the other hand, with respect to the $.$-product, $V$ is a
($\CC[T]$-)\emph{differential algebra} 
(i.e., an algebra with a derivation $T$)
with a unit $\vac$.

We summarize the above discussion in the following lemma.

\begin{lemma}\label{lem:c.0}
Giving a state--field correspondence on
a pointed vector space $(V,\vac)$ is equivalent to
providing $V$ with a structure of a $\CC[T]$-conformal algebra
and a structure of a\/ $\CC[T]$-differential algebra with a unit $\vac$.
\end{lemma}


Next, we translate the $n$-th product axioms in terms 
of the $\la$- and $.$-products.

\begin{lemma}\label{lem:c.1}
Let $(V, \vac)$ be a pointed vector space and let $Y$ be a
state--field correspondence. Fix $a,b,c \in V$.
Then the collection of $n$-th product identities
\eqref{eq:lln} for
$n\ge0$ implies
\begin{align}
\label{eq:c.8}
(a_\la b)_{\la+\mu} c &= a_\la(b_\mu c) - b_\mu(a_\la c) \, ,
\\
\label{eq:c.10}
a_\la(b.c) &= (a_\la b).c + b.(a_\la c)
  + \int_0^{\la} (a_\la b)_\mu c \, \di\mu \, .
\intertext{The $(-1)$-st product identity 
$Y(a_{(-1)} b, z) c = \bigl( Y(a,z)_{(-1)} Y(b,z) \bigr) c$
implies}
\label{eq:c.14}
(a.b)_\la c &= (e^{T \d_\la} a).(b_\la c) + (e^{T \d_\la} b).(a_\la c)
              + \int_0^\la b_\mu(a_{\la-\mu} c) \di\mu \, ,
\\
\label{eq:c.16}
(a.b).c - a.(b.c) &= 
   \Bigl( \int_0^T \di\la \, a \Bigr).(b_\la c)
 + \Bigl( \int_0^T \di\la \, b \Bigr).(a_\la c) \, .
\end{align}
\end{lemma}
\begin{proof}
The collection of $n$-th product identities \eqref{eq:lln} for $n\ge0$ 
is equivalent to:
\begin{equation}\label{eq:c.7}
Y(a_\la b, w)c = \Res_z e^{\la(z-w)} [Y(a,z), Y(b,w)] \, c
= e^{-\la w} [a_\la , Y(b,w)] \, c.
\end{equation}
Taking $\Res_w e^{(\la+\mu)w}$, we obtain \eqref{eq:c.8}.
Taking $\Res_w w^{-1}$, 
and using that 
$e^{-\la w} w^{-1} = w^{-1} + \tint_0^{-\la} e^{\mu w} \di\mu$,
we get:
\begin{equation*}\label{eq:c.9}
(a_\la b).c 
= a_\la(b.c) - b.(a_\la c)
  + \int_0^{-\la} [a_\la , b_\mu]c \, \di\mu \, .
\end{equation*}
This, together with \eqref{eq:c.8}, implies \eqref{eq:c.10}
(after the substitution $\mu' = \la+\mu$ in the integral).

Due to \eqref{eq:nop} and \eqref{eq:c.5}, the $(-1)$-st product
identity is equivalent to:
\begin{equation*}\label{eq:c.12}
Y(a.b, z)c = (e^{zT} a).(Y(b,z) c) 
             + Y(b,z) (a_{-\d_z} c) (z^{-1}) .
\end{equation*}
Taking $\Res_z e^{\la z}$ and using integration by parts, we get:
\begin{align*}
(a&.b)_\la c = \Res_z (e^{T \d_\la} e^{\la z} a).(Y(b,z) c)
            + \Res_z Y(b,z) (a_{\la-\d_z} c) (e^{\la z} z^{-1})
\\
&= (e^{T \d_\la} a).(b_\la c)
 + \Res_z Y(b,z) (a_{\la-\d_z} c) 
                 \Bigl( z^{-1} + \int_0^\la e^{\mu z} \di\mu \Bigr)
\\
&= (e^{T \d_\la} a).(b_\la c)
 + \Res_z \bigl( e^{\d_z \d_\la} Y(b,z) \bigr) (a_\la c) z^{-1}
 + \int_0^\la \Res_z Y(b,z) (a_{\la-\mu} c) e^{\mu z} \di\mu 
\, ,
\end{align*}
which implies \eqref{eq:c.14}, using translation invariance.

On the other hand, taking $\Res_z z^{-1}$
of the $(-1)$-st product identity, we get:
\begin{align*}
(a&.b).c = 
   \Res_z z^{-1} \bigl( Y(a,z)_+ Y(b,z)_+ c
 + Y(a,z)_+ Y(b,z)_- c
 + Y(b,z)_+ Y(a,z)_- c \bigr)
\\
&= a.(b.c)
 + \Res_z z^{-1} ((e^{zT}-1) a).(Y(b,z)_- c)
 + \Res_z z^{-1} ((e^{zT}-1) b).(Y(a,z)_- c)
\\
&= a.(b.c)
 + \Res_z \Bigl( \int_0^T e^{\la z} \di\la \, a \Bigr).(Y(b,z)_- c)
 + \Res_z \Bigl( \int_0^T e^{\la z} \di\la \, b \Bigr).(Y(a,z)_- c) \, ,
\end{align*}
which proves \eqref{eq:c.16}.
%
%
\end{proof}

\begin{remark}\label{rem:c.1}
Equation~\eqref{eq:c.7} is equivalent to the commutator formula
\begin{equation}\label{eq:c.11}
[a_{(m)}, b_{(n)}] \, c 
= \sum_{j=0}^\infty \binom{m}{j} (a_{(j)} b)_{(m+n-j)} c \, ,
\qquad m\in\ZZ_+ , \; n\in\ZZ \, .
\end{equation}
Equation~\eqref{eq:c.8} is equivalent to the same formula
for $m,n\in\ZZ_+$.
\end{remark}

Identity \eqref{eq:c.8} is called the 
(left) \emph{Jacobi identity} (cf.\ \cite[Sec.~2.7]{K}).
A conformal algebra satisfying this identity
for all $a,b,c \in V$
is called a (left) \emph{Leibniz conformal algebra}.
Equation~\eqref{eq:c.10} is known as
the ``non-commutative''  \emph{Wick formula} (cf.\ \cite[(3.3.12)]{K}),
while \eqref{eq:c.16} is called the \emph{quasi-associativity formula}
(cf.\ \cite[(4.8.5)]{K}).

Notice that the right-hand side of \eqref{eq:c.16} is symmetric
with respect to $a$ and $b$, hence \eqref{eq:c.16} implies
\begin{equation}\label{eq:c.16'}
a.(b.c) - b.(a.c) = (a.b - b.a).c \, .
\end{equation}
An algebra satisfying \eqref{eq:c.16'} 
for all $a,b,c \in V$
is called \emph{left-symmetric}.
For such an algebra $a.b - b.a$ is a Lie algebra bracket.

%
\begin{theorem}\label{thm:c.1}
Giving a strong field algebra structure on a pointed vector space $(V,\vac)$ 
is the same as providing $V$ with a structure of a 
Leibniz\/ $\CC[T]$-conformal algebra
and a structure of a\/ $\CC[T]$-differential algebra with a unit $\vac$,
satisfying \eqref{eq:c.10}--\eqref{eq:c.16}.
\end{theorem}
\begin{proof}
If $(V,\vac,Y)$ is a strong field algebra, then 
by the above discussion we can define 
a $\la$-product and a $.$-product on $V$ satisfying
all the requirements.

Conversely, given a $\la$-product and a $.$-product,
we define a state--field correspondence $Y$ by \eqref{eq:c.5}.
%
%
In the proof of \leref{lem:c.1},
we have seen that equations \eqref{eq:c.8}--\eqref{eq:c.16} are
equivalent to the identities
\begin{equation*}\label{eq:c.21}
\Res_z \bigl( Y(a_{(n)} b, z) - Y(a, z)_{(n)} Y(b, z) \bigr) F = 0,
\qquad a,b\in V, \; n\ge -1, \; F=e^{\la z} \text{ or } z^{-1}.
\end{equation*}

Using the translation invariance of $Y$ and integration by parts, 
we see that this identity holds also with $F$ replaced with $\d_z F$.
Hence it holds for all $F=z^k$, $k<0$. 
For $F=e^{\la z}$, taking coefficients at powers of $\la$ shows that
it is satisfied also for $F=z^k$, $k\ge0$.
This implies the $n$-th product axioms for $n\ge -1$. 

Replacing $a$ with $Ta$ and using translation invariance shows that,
for $n\ne0$, the $n$-th product axiom 
implies the $(n-1)$-st product axiom.
Therefore, the $n$-th product axioms hold for all $n\in\ZZ$,
and $(V,\vac,Y)$ is a strong field algebra.
\end{proof}

\begin{corollary}\label{cor:c.1}
Let $(V, \vac)$ be a pointed vector space, let $Y$ be a
state--field correspondence,
and define the $\la$-product and $.$-product on $V$ as above.
Let $U$ be a $T$-invariant subspace of $V$.

\alphaparenlist
\begin{enumerate}
\item 
If \eqref{eq:c.8} and \eqref{eq:c.10} hold for all $a,b\in U$, $c\in V$,
then $Y(a,z)_{(n)} Y(b,z) = Y(a_{(n)} b, z)$ 
for all $n\ge0$, $a,b\in U$.

\item 
If \eqref{eq:c.14} and \eqref{eq:c.16} hold for all $a,b\in U$, $c\in V$,
then $Y(a,z)_{(n)} Y(b,z) = Y(a_{(n)} b, z)$ 
for all $n<0$, $a,b\in U$.
\end{enumerate}
\end{corollary}

\section{Tensor Algebra over a Leibniz Conformal Algebra}
\label{sec:tens}

In this section, we are going to use the notation of the previous one.
Let $R$ be a Leibniz $\CC[T]$-conformal algebra (see \seref{sec:conf}),
and let $V = \tens(R) = \bigoplus_{m=0}^\infty R^{\tt m}$ 
be the tensor algebra over $R$ (viewed as a vector space over $\CC$). 
We extend the action of $T$ on $R$ to $V$ so that it is
a derivation of the tensor product. Let $\vac$ be the
element $1\in\CC \equiv R^{\tt 0} \subset V$.

\begin{theorem}\label{thm:t.1}
There exists a unique structure of 
a field algebra on $V = \tens(R)$ such that the restriction of the 
$\la$-product to $R \times R$ coincides with the $\la$-product of $R$,
the restriction of the $.$-product to $R \times V$ coincides with the
tensor product, and the $n$-th product axioms \eqref{eq:3.4} hold
for $a \in R$, $b \in V$.
\end{theorem}

For $a\in R$, $C\in V$, we set $a.C = a \tt C$. Next, we
define $a_\la C$ inductively using the Wick formula \eqref{eq:c.10}
and starting from $a_\la \vac = 0$:
%
\begin{equation}\label{eq:t.1}
a_\la(c \tt C) = (a_\la c) \tt C + c \tt (a_\la C)
  + \int_0^{\la} (a_\la c)_\mu C \, \di\mu \, ,
\qquad a,c\in R, \; C\in V.
\end{equation}
\begin{lemma}\label{lem:t.1}
Formula \eqref{eq:t.1} defines a representation of 
the Leibniz conformal algebra $R$ on $V$, 
whose restriction to $R$ coincides with the adjoint
representation.
\end{lemma}
\begin{proof}
For $a,b\in R$, $C\in V$, we have to check the translation invariance:
\begin{equation*}\label{eq:t.2}
(T a)_\la C = -\la \, a_\la C, \qquad
a_\la (T C) = (\la+T) (a_\la C),
\end{equation*}
and the Jacobi identity (cf.\ \eqref{eq:c.8}):
\begin{equation*}\label{eq:t.3}
(a_\la b)_{\la+\mu} C = a_\la(b_\mu C) - b_\mu(a_\la C) .
\end{equation*}
Note that these identities hold when $C\in R$ by definition.
Assuming that they hold for a fixed $C\in V$ and all
$a,b\in R$, we are going to prove them
for $a,b,c \tt C$ for any $a,b,c\in R$.

The translation invariance is immediate. The Jacobi identity reduces
to checking that
\begin{align*}
\int_0^{\la+\mu} & \bigl( (a_\la b)_{\la+\mu} c \bigr) {}_\nu C \,\di\nu
\\
&=
\int_0^\la \bigl( a_\la (b_\mu c) \bigr) {}_\nu C \,\di\nu
+
\int_0^\la (a_\la c)_\nu (b_\mu C) \,\di\nu
+
\int_0^\mu a_\la \bigl( (b_\mu c)_\nu C \bigr) \,\di\nu
\\
&-
\int_0^\mu \bigl( b_\mu (a_\la c) \bigr) {}_\nu C \,\di\nu
-
\int_0^\mu (b_\mu c)_\nu (a_\la C) \,\di\nu
-
\int_0^\la b_\mu \bigl( (a_\la c)_\nu C \bigr) \,\di\nu \, .
\end{align*}

By Jacobi identity, 
the third and fifth terms in the right-hand side combine to
$
\int_0^\mu \bigl( a_\la (b_\mu c) \bigr) {}_{\la+\nu} C \,\di\nu
=
\int_\la^{\la+\mu} \bigl( a_\la (b_\mu c) \bigr) {}_\nu C \,\di\nu
$,
which together with the first term gives
$
\int_0^{\la+\mu} \bigl( a_\la (b_\mu c) \bigr) {}_\nu C \,\di\nu
$.
Similarly, the sum of the other three terms is
\linebreak
$
-\int_0^{\la+\mu} \bigl( b_\mu (a_\la c) \bigr) {}_\nu C \,\di\nu
$.
Now the statement follows from the Jacobi identity for $a,b,c$.
\end{proof}

We have defined a $\la$-product $a_\la C$ and a $.$-product $a.C = a \tt C$
for $a\in R$, $C\in V$, satisfying the translation invariance, 
the Jacobi identity,  
and the Wick formula.  
By the results of \seref{sec:conf} (see \leref{lem:c.0} and \coref{cor:c.1}), 
this gives
a linear map from $R$ to $\glf(V)$, $a\mapsto Y(a,z)$,
such that $[T, Y(a,z)] = Y(Ta,z) = \d_z Y(a,z)$
and $Y(a,z)_{(n)} Y(b,z) = Y(a_{(n)} b, z)$ 
for all $n\ge0$, $a,b\in R$.

Next, we set $Y(\vac,z) = I_V$, and define fields $Y(A,z) \in \glf(V)$ 
inductively by the formula 
\begin{equation}\label{eq:t.5a}
Y(a \tt A,z) = \nop{Y(a,z) Y(A,z)}
\qquad \text{for \, $a\in R$, $A\in V$}. 
\end{equation}
It is easy to check that $Y$ is a
state--field correspondence.

%
%

%
\begin{lemma}\label{lem:t.2}
We have $Y(a_{(n)} B, z) = Y(a,z)_{(n)} Y(B,z)$
for all $a\in R$, $B\in V$, $n\in\ZZ$.
In particular, all pairs $(Y(a,z), Y(B,z))$ are weakly local
for $a\in R$, $B\in V$.
\end{lemma}
\begin{proof}
By definition, 
$Y(a_{(-1)} B, z) = Y(a \tt B, z) = \nop{Y(a,z) Y(B,z)} \, ,$
so the statement holds for $n=-1$ and hence for all $n<0$.
Thus, it suffices to show that 
\begin{equation*}\label{eq:t.6}
Y(a_\la B, z) = [Y(a,z)_\la Y(B,z)] 
= \Res_x e^{\la(x-z)} [Y(a,x), Y(B,z)], 
\qquad a\in R, \; B\in V.
\end{equation*}
We have already proved this when $B\in R$, while for $B=\vac$ both sides
are trivial.

Assuming the above identity is true for $B\in V$ and all $a\in R$, 
consider it for $b \tt B$
where $b\in R$. Using this assumption and \eqref{eq:t.1}, we find:
\begin{align*}
Y \bigl( a_\la(b \tt B), z \bigr) 
    = \bigl[ Y(a, z&)_\la Y(b, z) \bigr] {}_{(-1)} Y(B, z) 
    + Y(b, z)  {}_{(-1)} \bigl[ Y(a, z)_\la Y(B, z) \bigr]
\\
   &+ \int_0^{\la} 
          \bigl[ \bigl[ Y(a, z)_\la Y(b, z) \bigr] {}_\mu Y(B, z) \bigr] 
\, \di\mu \, .
\end{align*}
But the right-hand side is equal to 
$[ Y(a,z)_\la ( Y(b, z)_{(-1)} Y(B,z) ) ]$
by the Wick formula for arbitrary fields \cite[(3.3.12)]{K}.
This completes the proof.
\end{proof}

\begin{lemma}\label{lem:t.3}
All pairs $(Y(A,z), Y^\op(B,z))$ are local for $A,B \in V$.
\end{lemma}
\begin{proof}
It follows from \leref{lem:t.2} and \coref{cor:3.2} that
all pairs $(Y(a,z), Y^\op(B,z))$ are local for $a \in R$, $B \in V$.
By induction, assume that $(Y(A,z), Y^\op(B,z))$ is local, and
consider the pair $(Y(a \tt A, z), Y^\op(B,z))$ for $a \in R$.
Recall that by definition $Y(a \tt A, z) = Y(a,z)_{(-1)} Y(A,z)$.
Since, by \leref{lem:t.2}, the pair $(Y(a,z), Y(A,z))$
is weakly local, we can apply \leref{ldong}a to conclude that
the pair $(Y(a,z)_{(-1)} Y(A,z)$, $Y^\op(B,z))$ is local.
\end{proof}


Now \thref{th:3.1}a and \leref{lem:t.3} imply that
the so defined $(V,\vac,Y)$ is a field algebra.
Uniqueness is clear by construction.
This completes the proof of \thref{thm:t.1}.

The field algebra $V=\tens(R)$
has the following universality property.
Let $W$ be a strong field algebra and let $f \colon R\to W$
be a homomorphism of conformal algebras
(i.e., $f(a_{(n)} b) = f(a)_{(n)} f(b)$ and $f(Ta) = T f(a)$
for all $a,b \in R$, $n \in \ZZ_+$).
Then there is a unique homomorphism of field algebras
$\wti f \colon V\to W$ such that $f = \wti f \circ i$,
where $i$ is the embedding of $R$ in $V$. 

\begin{remark}\label{rem:t.1}
Although all pairs $(Y(a,z), Y(B,z))$ are 
weakly local for $a \in R$, $B \in V$ (see \leref{lem:t.2}),
in general $(V,\vac,Y)$ is not a strong field algebra.
This follows from the results of \seref{sec:sfv}.
\end{remark}

\begin{remark}\label{rem:g.4}
A conformal algebra $R$ is called \emph{graded}
if there is a diagonalizable operator $H \in \End R$ satisfying 
$[H,T]=T$ and \eqref{eq:g.2} for $n\in\ZZ_+$ and homogeneous $a \in R$.
Then the tensor algebra $\tens(R)$ over a Leibniz conformal algebra $R$ 
is graded if $R$ is graded.
Moreover, $\tens(R)$ is $\ZZ_+$-graded if $R$ is $\ZZ_+$-graded.
\end{remark}

\section{Field Algebras and Vertex Algebras}
\label{sec:vert}

Recall that a \emph{vertex algebra} $(V,\vac,Y)$
is a pointed vector space $(V,\vac)$ together with
a local state--field correspondence $Y$ (see e.g.~\cite{K}).
In particular, any vertex algebra is a strong field algebra
with $Y=Y^\op$ (see \thref{th:3.1}).

\begin{example}
  \label{ex:3.1}

The state--field correspondence given by Example~\ref{ex:1.1} is
a vertex algebra iff the algebra $V$ is commutative and
associative.

\end{example}

\begin{example}
  \label{ex:3.2}

Let $V$ be a subspace of $\glf(U)$ which contains $I_U$,
is $\partial_x$-invariant and is closed under all $n$-th
products.  Then the state--field correspondence defined in
\exref{ex:1.2} gives $V $ a structure of a vertex algebra
iff $V$ is a local collection.  This follows from 
\cite[Proposition~3.2]{K}.

\end{example}

\begin{theorem}\label{thm:vert}


  A vertex algebra is the same as a field algebra $(V,\vac,Y)$ 
  for which $Y=Y^\op$.


\end{theorem}
\begin{proof}
Let $(V,\vac,Y)$ be a field algebra with $Y=Y^\op$.
Recall that (see \prref{prop:ax}) 
the associativity relation \eqref{eq:3.9} is equivalent to: 
\begin{equation*}
(z-w)^N [Y(a,z), Y^\op(c,w)] \, b = 0.
\end{equation*}
By \leref{lem:2.2a}, this relation also holds after replacing
$b$ with $Tb$.
Replacing $b$ with $e^{uT}b$ in \eqref{eq:3.9}
and using \prref{prop:1.1}b, we get:
\begin{equation}\label{eq:asu1}
\begin{split}
  (z-w&)^N i_{z,u} i_{w,u} Y(Y(a,z-u) b,u-w) c
\\ 
&= (z-w)^N i_{z,w} i_{w,u} Y(a,z-w) Y(b,u-w) c \, .
\end{split}
\end{equation}

There exists $P \in \ZZ_+$ such that $b_{(k)} c=0$ for $k \geq P$.
Hence the right-hand side of \eqref{eq:asu1} multiplied by
$(u-w)^P$ contains only non-negative powers of $u-w$, hence only
non-negative powers of $w$.  Therefore, multiplying \eqref{eq:asu1}
by $(u-w)^P$ and putting $w=0$ gives:
\begin{equation}\label{eq:asu2}
\begin{split}
  Y(a,z&) Y(b,u) c
\\ 
  &= z^{-N} u^{-P} (z-w)^N (u-w)^P
     i_{z,u} i_{w,u} Y(Y(a,z-u) b, u-w) c \big|_{w=0} \, .
\end{split}
\end{equation}

Using again \prref{prop:1.1}b and $Y=Y^\op$,
we compute:
\begin{equation*}
i_{z,u} Y(a, z-u) b = Y(e^{-uT}a, z) b 
= e^{zT} Y(b, -z) e^{-uT}a
= i_{z,u} e^{(z-u)T} Y(b, u-z) a \, .
\end{equation*}
Therefore,
\begin{equation}\label{eq:asu3}
i_{z,u} i_{w,u} Y(Y(a,z-u) b, u-w)
= i_{z,u} i_{w,z} Y(Y(b,u-z) a, z-w) \, .
\end{equation}

We may take $N=P\gg0$ in \eqref{eq:asu2}.
Comparing \eqref{eq:asu2} and \eqref{eq:asu3}, we see that 
$Y(b,u) Y(a,z) c$ is given by the right-hand side of \eqref{eq:asu2}
where $i_{z,u}$ is replaced with $i_{u,z}$. Hence, if $K\in\ZZ_+$
is such that $a_{(k)} b = 0$ for $k\ge K$, we have
\begin{equation*}
(z-u)^K Y(a,z) Y(b,u) c = (z-u)^K Y(b,u) Y(a,z) c \, ,
\end{equation*}
which is 
locality of $Y(a,z)$ and $Y(b,z)$.
This proves that $(V,\vac,Y)$ is a vertex algebra.
\end{proof}

\begin{corollary}\label{cor:3.1}
A state--field correspondence is local iff it is local on every vector.


\end{corollary}

\begin{remark}\label{rem:3.2}
It follows from \eqref{eq:asu2}, \eqref{eq:asu3} that
for $N\gg0$ and any $n\in\ZZ$:
\begin{equation}\label{eq:asu4a}
\begin{split}
  Y(a,z) Y(b,u) \, i_{z,u} &(z-u)^n 
- Y(b,u) Y(a,z) \, i_{u,z} (z-u)^n 
\\ 
= z^{-N} u^{-N} &(z-w)^N (u-w)^N i_{w,u} 
   (i_{z,u} - i_{u,z}) (z-u)^n \,
\\
&\times 
  Y(Y(a,z-u) b, u-w) c \big|_{w=0} \, .
\end{split}
\end{equation}
Using \eqref{eq:1.2}, we find
\begin{equation*}
(i_{z,u} - i_{u,z}) (z-u)^n \, Y(a,z-u) b
= \sum_{m=0}^\infty  
  ( a_{(m+n)} b ) \, \d_u^m \de(z-u) / m! \, .
\end{equation*}
Note that this sum is finite, because $a_{(m)} b = 0$ for
$m\gg0$. Therefore, for large enough $N$, we can put
$w=0$ in \eqref{eq:asu4a}, and obtain:
\begin{equation}\label{eq:asu4}
\begin{split}
  Y(a,z&) Y(b,u) \, i_{z,u} (z-u)^n 
- Y(b,u) Y(a,z) \, i_{u,z} (z-u)^n 
\\ 
&= \sum_{m=0}^\infty  
  Y(a_{(m+n)} b, u) \, \d_u^m \de(z-u) / m! \, .
\end{split}
\end{equation}
The collection of these identities for $n\in\ZZ$ is 
equivalent to the Borcherds identity \cite[(4.8.1)]{K}.
In particular, taking $\Res_z$ of both sides of \eqref{eq:asu4},
we obtain the $n$-th product identity \eqref{eq:3.4}.
\end{remark}

\begin{remark}\label{rem:3.3}
Note that formula \eqref{eq:asu1} holds for any field algebra $(V,\vac,Y)$.
Let $K \in \ZZ_+$ be such that $a_{(k)} b=0$ for $k \geq K$.
After multiplication by $(z-u)^K$, the left-hand side of \eqref{eq:asu1}
contains only non-negative powers of $z-u$, hence only
non-negative powers of $z$. Putting $z=0$
and replacing $u$ with $-u$, we obtain:
\begin{equation*}
\begin{split}
i_{w,u} & Y(Y(a,u) b,-u-w) c
\\ 
&= (-w)^{-N} u^{-K} (z-w)^N (z+u)^K
   i_{z,w} i_{w,u} Y(a,z-w) Y(b,-u-w) c \big|_{z=0} \, .
\end{split}
\end{equation*}
Since in both sides $-u-w$ is expanded in non-negative powers of $u$, 
we can replace $w$ with $-u-w$ and use Taylor's formula to get:
\begin{equation}\label{eq:asu6}
\begin{split}
Y(Y(a,u) b,w) &c
\\ 
= i_{z,u+w} & i_{w,u} (u+w)^{-N} u^{-K} (z+u+w)^N (z+u)^K
\\
&\times Y(a,z+(u+w)) Y(b,w) c \big|_{z=0} \, .
\end{split}
\end{equation}
\end{remark}

\begin{remark}\label{rem:6.7a}

Let $V$ be a field algebra. Given $a,b \in V$, denote by $ab$
the $\CC$-span of all elements $a_{(n)} b$ ($n \in \ZZ$).
More generally, if $A$ and $B$ are two subspaces of $V$,
we define $AB$ as the $\CC$-span of all elements $a_{(n)} b$ 
($n \in \ZZ$, $a \in A$, $b \in B$).
Formulas \eqref{eq:asu2} and \eqref{eq:asu6} 
imply associativity of this product:
\begin{equation}\label{eq:6.5}
(ab)c = a(bc) \, , \qquad a,b,c \in V \, .
\end{equation}
This property allows one to apply many arguments of the theory of
associative algebras to field algebras. 
Notice that \eqref{eq:6.5} holds also for $a,b \in V$ and $c \in M$,
where $M$ is a module over the field algebra $V$ (see \deref{def:5.1} below).

\end{remark}

\begin{remark}\label{rem:6.7}
Letting $a \bullet b = \CC[T](ab)$, we still have associativity of the product
$\bullet$ for an arbitrary field algebra $V$. In addition, 
$\vac \bullet a = a \bullet \vac = \CC[T]a$ for all $a \in V$.
If $V$ is a vertex algebra, the product $\bullet$ is also commutative.
These remarks allow one to use arguments of commutative algebra
to study vertex algebras
(this will be pursued in a future publication).
\end{remark}

For the opposite state--field correspondence $Y^\op$,
we can consider
\begin{equation*}\label{eq:c.22}
a_\la^\op b = \Res_z e^{\la z} Y^\op(a,z)b,
\qquad
a \overset{\op}. b = \Res_z z^{-1} Y^\op(a,z)b,
\end{equation*}
and express these operations
in terms of the $\la$-product and $.$-product
defined for $Y$. We have:
\begin{align}\label{eq:c.23}
a_\la^\op b &= \Res_z e^{\la z} e^{zT} Y(b,-z)a
   = -\Res_z e^{-(\la+T)z} Y(b,z)a
   = -b_{-\la-T} a,
\intertext{and, using 
$z^{-1} e^{zT} = z^{-1} - \int_0^{-T} e^{-\la z} \, \di\la\,$,}
\label{eq:c.24}
a \overset{\op}. b &= \Res_z z^{-1} e^{zT} Y(b,-z)a
= b.a + \int_0^{-T} b_\la a \, \di\la \, .
\end{align}
%
%

If $(V,\vac,Y)$ is a vertex algebra, then the
$\la$-product satisfies the \emph{skewsymmetry relation}
$a_\la b = -b_{-\la-T} a$ for $a,b\in V$.
A conformal algebra satisfying the Jacobi identity and
the skewsymmetry relation is called a \emph{Lie conformal algebra}
\cite[Sec.~2.7]{K}.

\begin{theorem}\label{thm:c.2}
Giving a vertex algebra structure on a pointed vector space $(V,\vac)$ 
is the same as providing $V$ with the structures of a 
Lie\/ $\CC[T]$-conformal algebra
and a\/ left-symmetric 
$\CC[T]$-differential algebra with a unit $\vac$
{\rm(}see \eqref{eq:c.16'}{}{\rm)},
satisfying the Wick formula \eqref{eq:c.10} and
\begin{equation}\label{eq:c.25}
a.b - b.a = \int_{-T}^0 a_\la b \, \di\la \, ,
\qquad a,b \in V.
\end{equation}
\end{theorem}
\begin{proof}
By \leref{lem:c.0}, these data define a state--field correspondence $Y$.
By \eqref{eq:c.23} and \eqref{eq:c.24}, we have $Y=Y^\op$;
therefore, we only need to prove that $V$ is a field algebra
(see \thref{thm:vert}).
Due to \thref{thm:c.1}, it remains to show that
\eqref{eq:c.14} and \eqref{eq:c.16} follow from 
equations~\eqref{eq:c.8}, \eqref{eq:c.10}, \eqref{eq:c.16'},
\eqref{eq:c.25} and the skewsymmetry of the $\la$-product.

First, using \eqref{eq:c.25}, rewrite the Wick formula \eqref{eq:c.10} as
\begin{equation*}\label{eq:c.26}
a_\la(b.c) = c.(a_\la b) + b.(a_\la c)
  + \int_{-T}^{\la} (a_\la b)_\mu c \, \di\mu \, .
\end{equation*}
Then replace $\la$ with $-\la$ in this equation and apply $e^{T \d_\la}$
to both sides. Using Taylor's formula and the fact 
that $T$ is a derivation, we obtain:
\begin{equation*}\label{eq:c.27}
a_{-\la-T}(b.c) = (e^{T \d_\la} c).(a_{-\la-T} b) 
  + (e^{T \d_\la} b).(a_{-\la-T} c)
  + \int_{-T}^{-\la-T} e^{T \d_\la} 
             \bigl( (a_{-\la} b)_\mu c \bigr) \, \di\mu \, .
\end{equation*}
Using \eqref{eq:c.4} and the skewsymmetry of the $\la$-product,
the last term can be rewritten as follows:
\begin{align*}
\int_{-T}^{-\la-T} e^{T \d_\la} 
             \bigl( (a_{-\la} b)_\mu c \bigr) \, \di\mu 
&= - \int_{-T}^{-\la-T} e^{T \d_\la}
             \bigl( c_{-\mu-T} (a_{-\la} b) \bigr) \, \di\mu
\\
&= \int_0^\la e^{T \d_\la} \bigl( c_\mu (a_{-\la} b) \bigr) \, \di\mu
= \int_0^\la c_\mu (a_{-\la+\mu-T} b) \, \di\mu \, .
\end{align*}
Applying the skewsymmetry relation again, we obtain \eqref{eq:c.14}.

In order to prove \eqref{eq:c.16}, we manipulate its left-hand side
using twice \eqref{eq:c.25}:
\begin{equation*}\label{eq:c.29}
(a.b).c - a.(b.c) = c.(a.b) - a.(c.b)
   + \int_{-T}^0 \di\la\, (a.b)_\la c
   - a. \Bigl( \int_{-T}^0 \di\la\, b_\la c \Bigr).
\end{equation*}
Due to \eqref{eq:c.16'} and \eqref{eq:c.25}, 
the first two terms in the right-hand side of this equation give:
\begin{align*}
c.(a.b) - a.(c.b) &= -(a.c-c.a).b 
   = - \Bigl( \int_{-T}^0 \di\la\, a_\la c \Bigr) .b
\\  
  &= - b. \Bigl( \int_{-T}^0 \di\la\, a_\la c \Bigr)
    +  \int_{-T}^0 \di\mu\, b_\mu\Bigl( \int_{-T}^0 \di\la\, a_\la c \Bigr) .
\end{align*}
The double integral is equal to:
\begin{equation}\label{eq:c.31}
\begin{split}
\int_{-T}^0 & \di\mu\, b_\mu\Bigl( \int_{-T}^0 \di\la\, a_\la c \Bigr) 
   = \int_{-T}^0 \di\mu\, \int_{-T-\mu}^0 \di\la\, b_\mu(a_\la c)
\\
  &= \int_{-T}^0 \di\mu\, \int_{-T}^\mu \di\la\, b_\mu(a_{\la-\mu} c)
   = \int_{-T}^0 \di\la\, \int_\la^0 \di\mu\, b_\mu(a_{\la-\mu} c) .
\end{split}
\end{equation}
Therefore, applying \eqref{eq:c.14}, which we already proved, we obtain:
\begin{align*}
(a.b).c - a.(b.c) &= 
    \int_{-T}^0 \di\la\, (e^{T \d_\la} a).(b_\la c)
  - a. \Bigl( \int_{-T}^0 \di\la\, b_\la c \Bigr)
\\
  &+ \int_{-T}^0 \di\la\, (e^{T \d_\la} b).(a_\la c)
  - b. \Bigl( \int_{-T}^0 \di\la\, a_\la c \Bigr) .
\end{align*}
{}From here, it is easy to deduce \eqref{eq:c.16}.

This completes the proof of the theorem.
\end{proof}

\begin{proposition}\label{prop:c.2}
If $R$ is a Leibniz\/ $\CC[T]$-conformal algebra with a
$\la$-product $a_\la b$, then the bracket
\begin{equation}\label{eq:c.33}
[a,b] = \int_{-T}^0  a_\la b \, \di\la \, ,
\qquad a,b \in R
\end{equation}
defines the structure of a Leibniz\/ $\CC[T]$-differential algebra on $R$.
If $R$ is a Lie conformal algebra, then this is a Lie bracket.
\end{proposition}
\begin{proof}
It is clear that $T$ is a derivation of the bracket \eqref{eq:c.33},
because it is a derivation of the $\la$-product $a_\la b$.
{}From relation \eqref{eq:c.31}, we have:
\begin{align*}
[b,[a,c]]
 & 
  = \int_{-T}^0 \di\la\, \int_\la^0 \di\mu\, b_\mu(a_{\la-\mu} c) 
  = \int_{-T}^0 \di\la\, \int_\la^0 \di\mu\, b_{\la-\mu}(a_\mu c) \, ,
\intertext{and similarly,}
[a,[b,c]]
 &= \int_{-T}^0 \di\la\, \int_\la^0 \di\mu\, a_\mu(b_{\la-\mu} c) \, .
%
\intertext{On the other hand,}
[[a,b],c]
 &= \int_{-T}^0 \di\la\, \Bigl( 
                \int_{-T}^0 \di\mu\, a_\mu b \Bigr) {}_\la c
  = \int_{-T}^0 \di\la\, \int_\la^0 \di\mu\, (a_\mu b)_\la c \, .
\end{align*}
Therefore, the Jacobi identity for the bracket \eqref{eq:c.33}
follows from the Jacobi identity \eqref{eq:c.8} for the
$\la$-product.
Finally, it is easy to see that $[a,b] = -[b,a]$ if
$a_\la b = -b_{-\la-T} a$.
\end{proof}

\begin{remark}\label{rem:v.1}
Let $R$ be a Lie\/ $\CC[T]$-conformal algebra with a
$\la$-product $a_\la b = \sum_{m\ge0} \la^m \, a_{(m)} b / m! \,$.  
To it one can associate a Lie algebra $\Lie R$ as follows
(see \cite[Sec.~2.7]{K}). As a vector space, 
$\Lie R = R[t,t^{-1}] / (T+\d_t) R[t,t^{-1}]$.
The Lie bracket in $\Lie R$ is given by the formula
(cf.\ \eqref{eq:c.11}):
\begin{equation*}
[a_{[m]}, b_{[n]}] 
= \sum_{j=0}^\infty \binom{m}{j} (a_{(j)} b)_{[m+n-j]} \, ,
\qquad a,b \in R \, , \;\; m,n\in\ZZ \, ,
\end{equation*}
where $a_{[m]}$ is the image of $a t^m$ in $\Lie R$.
Note that $T$ acts on $\Lie R$ as a derivation: 
$T(a_{[m]}) = -m \, a_{[m-1]}$.

Let $(\Lie R)_-$ (respectively, $(\Lie R)_+$) be the $\CC$-span of all
$a_{[m]}$ for $a \in R$, $m \ge 0$ (respectively, $m < 0$).
These are subalgebras of $\Lie R$.
Let $R_{\mathrm{Lie}}$ be $R$ considered as a Lie algebra 
with respect to the bracket \eqref{eq:c.33}.
Then the map $a \mapsto a_{[-1]}$ is an isomorphism of
$\CC[T]$-differential Lie algebras $R_{\mathrm{Lie}} \isoto (\Lie R)_+$.
\end{remark}

\begin{theorem}\label{thm:v.2}
Let $R$ be a Lie\/ $\CC[T]$-conformal algebra with a
$\la$-product $a_\la b$. Let $R_{\mathrm{Lie}}$ be $R$ 
considered as a Lie algebra 
with respect to the bracket \eqref{eq:c.33}, and let 
$V=U(R_{\mathrm{Lie}})$ be its
universal enveloping algebra. Then there exists a unique structure of 
a vertex algebra on $V$ such that the restriction of the 
$\la$-product to $R_{\mathrm{Lie}} \times R_{\mathrm{Lie}}$ 
coincides with the $\la$-product of $R$
and the restriction of the $.$-product to 
$R_{\mathrm{Lie}} \times V$ coincides with the
product in $U(R_{\mathrm{Lie}})$.
\end{theorem}
\begin{proof}
We have $R_{\mathrm{Lie}} \simeq (\Lie R)_+$  (see \reref{rem:v.1}).
By construction, 
\begin{equation*}
V \simeq U((\Lie R)_+) \simeq \Ind^{\Lie R}_{(\Lie R)_-} \CC
\end{equation*}
(where $(\Lie R)_-$ acts trivially on $\CC$)
coincides with the universal vertex algebra associated to $\Lie R$
(see \cite[Sec.~4.7]{K}).
\end{proof}

The vertex algebra $V$ from \thref{thm:v.2}
has the following universality property.
Let $W$ be another vertex algebra and let $f \colon R\to W$
be a homomorphism of conformal algebras
(i.e., $f(a_{(n)} b) = f(a)_{(n)} f(b)$ and $f(Ta) = T f(a)$
for all $a,b \in R$, $n \in \ZZ_+$).
Then there is a unique homomorphism of vertex algebras
$\wti f \colon V\to W$ such that $f = \wti f \circ i$,
where $i$ is the embedding of $R$ in $V$. 
%
For this reason, $V$ is called the 
\emph{universal vertex algebra} associated to
$R$ and is denoted by $\univ(R)$.

\begin{remark}\label{rem:v.5}
$\univ(R)$ is the quotient of the field algebra $\tens(R)$,
constructed in \seref{sec:tens}, by the two-sided ideal generated by 
$a \tt b - b \tt a - \int_{-T}^0 a_\la b \, \di\la$
$(a,b \in R)$.
A different proof of \thref{thm:v.2} was given in \cite{GMS}.
\end{remark}

\section{Strong Field Algebras and Vertex Algebras}
\label{sec:sfv}

Let $(V,\vac,Y)$ be a strong field algebra.
We are going to use the $\la$-product and $.$-product on $V$
defined in \seref{sec:conf}.


Fix four elements $a,a',b,c \in V$.
Using formulas \eqref{eq:c.10}, \eqref{eq:c.14},
it is straightforward to compute:
\begin{align*}
(a.a')_\la (b_\mu c) 
 &= (e^{T \d_\la} a).(a'_\la (b_\mu c)) 
 + (e^{T \d_\la} a').(a_\la (b_\mu c))
 + \int_0^\la \di\nu\, a'_\nu (a_{\la-\nu} (b_\mu c))  \, ,
\\
\\
b_\mu ((a.a')_\la c)
 &= (e^{T \d_\la} a).(b_\mu (a'_\la c))
 + (e^{T \d_\la} (b_\mu a)).(a'_{\la+\mu} c)
\\
 + &\int_0^\mu \di\nu\, (b_\mu a)_\nu (a'_{\la+\mu-\nu} c)
 + (a \otto a')
 + \int_0^\la \di\nu\, b_\mu (a'_\nu (a_{\la-\nu} c))  \, ,
\intertext{and}
((a.a')_\la b)_{\la+\mu} c 
 &= (e^{T \d_\la} a).((a'_\la b)_{\la+\mu} c)
 + (e^{T \d_\la} (a'_{-\mu-T} b)).(a_{\la+\mu} c)
\\
 + &\int_0^{\la+\mu} \di\nu\, (a'_{\nu-\mu} b)_\nu (a_{\la+\mu-\nu} c)
 + (a \otto a')
 + \int_0^\la \di\nu\, (a'_\nu (a_{\la-\nu} b))_{\la+\mu} c \, ,
\end{align*}
where $(a \otto a')$ denotes the sum of the first three terms
with $a$ and $a'$ interchanged. Using Jacobi identity~\eqref{eq:c.8} 
for elements of the form $u$, $u_{(m)} v$ 
for $u,v \in\{ a,a',b,c \}$, $m\in\ZZ_+$, 
after a long calculation one obtains:
\begin{equation}\label{eq:sf.1}
\begin{split}
((a.a')_\la b)_{\la+\mu} c  &- (a.a')_\la (b_\mu c) + b_\mu ((a.a')_\la c)
\\
&= (e^{T \d_\la} (a_{-\mu-T} b + b_\mu a)).(a'_{\la+\mu} c)
 + (a \otto a') \, .
\end{split}
\end{equation}

This implies the following lemma.

\begin{lemma}\label{lem:sf.1}
Let $V$ be a vector space equipped with a structure of
a\/ $\CC[T]$-conformal algebra
and a structure of a\/ $\CC[T]$-differential algebra.
Let $a,a',b,c \in V$, and assume that formulas
\eqref{eq:c.8}--\eqref{eq:c.14} hold for all 
elements of the form $u$, $u_{(m)} v$ 
$(u,v \in\{ a,a',b,c \}, \, m\in\ZZ_+)$.
Then Jacobi identity~\eqref{eq:c.8} holds for the triple $(a.a',b,c)$
if and only if
\begin{equation}\label{eq:sf.2}
(a_{-\mu-T} b + b_\mu a).(a'_\la c) + (a \otto a') = 0 \, .
\end{equation}
\end{lemma}
\begin{proof}
By the above calculation, Jacobi identity for $(a.a',b,c)$
is equivalent to the vanishing of the right-hand side of \eqref{eq:sf.1}.
%
Applying to it $e^{-T \d_\la}$ and using Taylor's formula, 
we obtain:
\begin{equation*}
(a_{-\mu-T} b + b_\mu a)).(a'_{\la+\mu-T} c) 
+ (a \otto a')
= 0 \, ,
\end{equation*}
which is equivalent to \eqref{eq:sf.2}.
\end{proof}

Applying first \eqref{eq:c.14} and then \eqref{eq:c.10}, we find:
\begin{align*}
(a.a')_\la (b.c) 
 &= (e^{T \d_\la} a).\Bigl( (a'_\la b).c + b.(a'_\la c) 
           + \int_0^\la \di\mu\, (a'_\la b)_\mu c \Bigr)
  + \int_0^\la \di\mu\, (a_\mu b).(a'_{\la-\mu} c)
\\
 &+ (a \otto a')
  + \int_0^\la \di\mu\, 
  \Bigl\{ (a'_\mu (a_{\la-\mu} b)).c + b.(a'_\mu (a_{\la-\mu} c))   \Bigr\}
\\
 &+ \int_0^\la \di\mu\, \int_0^\mu \di\nu\,
  \Bigl\{ (a'_\mu (a_{\la-\mu} b))_\nu c 
          + (a'_\mu b)_\nu (a_{\la-\mu} c)   \Bigr\}
\\  
 &+ \int_0^\la \di\mu\, \int_0^{\la-\mu} \di\nu\,
  a'_\mu ((a_{\la-\mu} b)_\nu c) \, ,
\end{align*}
where $(a \otto a')$ denotes the sum of the first two terms
with $a$ and $a'$ interchanged. We also have:
\begin{align*}
((a.a')_\la b).c 
 &= \Bigl( (e^{T \d_\la} a).(a'_\la b) + (e^{T \d_\la} a').(a_\la b)
           + \int_0^\la \di\mu\, a'_\mu (a_{\la-\mu} b) \Bigr).c
\intertext{and}
b.((a.a')_\la c)
 &= b.\Bigl( (e^{T \d_\la} a).(a'_\la c) + (e^{T \d_\la} a').(a_\la c)
           + \int_0^\la \di\mu\, a'_\mu (a_{\la-\mu} c) \Bigr) \, .
\end{align*}
Using Jacobi identity~\eqref{eq:c.8}
and quasi-associativity~\eqref{eq:c.16}
for elements of the form $T^n u$, $u_{(m)} v$ 
for $u,v \in\{ a,a',b,c \}$, $m,n\in\ZZ_+$, 
after a long calculation one can show that:
\begin{equation}\label{eq:sf.3}
\begin{split}
(a&.a')_\la (b.c) - ((a.a')_\la b).c - b.((a.a')_\la c)
  - \int_0^\la \di\mu\, ((a.a')_\la b)_\mu c
\\
 &= \Bigl( (e^{T \d_\la} a).b - b.(e^{T \d_\la} a) 
       - \int_{-T}^0 \di\mu\, (e^{T \d_\la} a)_\mu b \Bigr).(a'_\la c) 
  + (a \otto a') \, .
\end{split}
\end{equation}

This implies:

\begin{lemma}\label{lem:sf.2}
Let $V$ be a vector space equipped with a structure of
a\/ $\CC[T]$-conformal algebra
and a structure of a\/ $\CC[T]$-differential algebra.
Let $a,a',b,c \in V$, and assume that formulas
\eqref{eq:c.8}--\eqref{eq:c.16} hold for all 
elements of the form $T^n u$, $u_{(m)} v$ 
$(u,v \in\{ a,a',b,c \}$, $m,n\in\ZZ_+)$.
Then the Wick formula \eqref{eq:c.10} holds for $(a.a',b,c)$
if and only if
\begin{equation}\label{eq:sf.4}
\Bigl( a.b - b.a - \int_{-T}^0 \di\mu\, a_\mu b \Bigr).(a'_\la c) 
  + (a \otto a') = 0 \, .
\end{equation}
\end{lemma}
\begin{proof}
By the above calculation, Wick formula for $(a.a',b,c)$
is equivalent to the vanishing of the right-hand side of \eqref{eq:sf.3}.
Note that \eqref{eq:c.10} holds for $(a,b,c)$ if and only if
it holds for $(Ta,b,c)$. Replacing in \eqref{eq:sf.3} $a$ and $a'$ by 
$e^{-T \d_\la} a$ and $e^{-T \d_\la} a'$, respectively, we obtain
\eqref{eq:sf.4}.
\end{proof}

The above two lemmas, together with \thref{thm:c.1} and
relations \eqref{eq:c.23}, \eqref{eq:c.24},
immediately imply the following:

\begin{corollary}\label{cor:sf.1}
Let $(V,\vac,Y)$ be a strong field algebra. Then all
$a,a',b,c \in V$ satisfy
\begin{equation}\label{eq:sf.5}
\Bigl( Y(a,z)b - Y^\op(a,z)b \Bigr).(a'_\la c) 
  + (a \otto a') = 0 \, .
\end{equation}
\end{corollary}

{}From this, we deduce the following surprising result.

\begin{theorem}\label{thm:sf.1}
Let $(V,\vac,Y)$ be a strong field algebra such that
$\vac = a_{(n)} c$ for some $a,c \in V$, $n \in \ZZ_+$.
Then $V$ is a vertex algebra.

In particular, if $V$ contains a Virasoro element, i.e.,
an element $\ell$ such that
$\ell_\la \ell = (T+2\la) \ell + \al \la^3 \vac$,
with $\al \in\CC\setminus\{0\}$, then $V$ is a vertex algebra.
\end{theorem}
\begin{proof}
Taking $a=a'$ and coefficient in front of $\la^n$ in \eqref{eq:sf.5},
we obtain $(Y(a,z)b - Y^\op(a,z)b).\vac = 0$, hence 
$Y(a,z)=Y^\op(a,z)$. Then \eqref{eq:sf.5} implies 
$Y(a',z)=Y^\op(a',z)$ for all $a' \in V$.
By \thref{thm:vert}, $V$ is a vertex algebra.
\end{proof}

\section{Modules over Field Algebras}
\label{sec:mod}

\begin{definition}\label{def:5.1}
A \emph{module} over a field algebra $(V,\vac,Y)$ is a 
vector space $M$ equipped with an operator $T^M \in \End M$
and a linear map 
\begin{equation*}
V \to \glf(M), \quad a \mapsto Y^M (a,z) = \sum_{n
  \in \ZZ} a^M_{(n)} z^{-n-1}, 
\end{equation*}
satisfying:
\begin{list}{}{}
\item (\emph{vacuum axiom}) \;
$Y^M (\vac ,z) =I_M$,

\item (\emph{translation invariance}) \;
$[T^M, Y^M(a,z)] = \partial_z Y^M(a,z)$,

\item (\emph{associativity axiom}) \;
  $(z-w)^N Y^M (Y(a,z)b,-w) v$ 
\\
  \hspace*{120pt}
  $= (z-w)^N i_{z,w} Y^M (a,z-w) Y^M (b,-w) v$
\end{list}
for $a,b \in V$, $v\in M$, $N \gg 0$.

When $(V,\vac,Y)$ is a strong field algebra, in the above definition
we replace the associativity axiom with the following
\begin{list}{}{}
\item (\emph{$n$-th product axiom}) \;
$Y^M(a_{(n)} b, z) = Y^M(a,z)_{(n)} Y^M(b,z)$,
\quad $n\in\ZZ$.
\end{list}
\end{definition}

\begin{remark}\label{rem:5.3}
What we introduced in \deref{def:5.1} are left modules over a field algebra.
We can define a \emph{right module} over a field algebra $(V,\vac,Y)$
as a left module over the opposite field algebra $(V,\vac,Y^\op)$.
Any field algebra is both a left and a right module over itself.
\end{remark}

Note that, as in \prref{prop:ax}, the $n$-th product axiom implies
the associativity axiom.
Taking $b=\vac$ in the associativity axiom
and using \prref{prop:1.1}b, we deduce
\begin{equation}\label{eq:5.2}
  Y^M (Ta,z) = \partial_z Y^M (a,z) \, .
\end{equation}
%

\begin{remark}\label{rem:5.2}
Let $(V,\vac,Y)$ be a vertex algebra, and let $M$ be a module
over $V$ considered as a field algebra. Then as in \seref{sec:vert}
one can show that (see \reref{rem:3.2}):
\begin{equation}\label{eq:asu5}
\begin{split}
  Y^M(a,z&) Y^M(b,u) \, i_{z,u} (z-u)^n 
- Y^M(b,u) Y^M(a,z) \, i_{u,z} (z-u)^n 
\\ 
&= \sum_{m=0}^\infty  
  Y^M(a_{(m+n)} b, u) \, \d_u^m \de(z-u) / m! \, ,
\qquad a,b \in V, \; n \in\ZZ \, .
\end{split}
\end{equation}
Therefore, $M$ is also a module over $V$ considered as a vertex algebra. 
\end{remark}

Let $M$ be a module over a field algebra $(V,\vac,Y)$.
%
%
%
%
Using \eqref{eq:5.2},
the same argument as in \seref{sec:vert} gives the following 
generalization of formula~\eqref{eq:asu2}:
%
%
%
\begin{equation}\label{eq:5.4}
\begin{split}
  Y^M &(a,z) Y^M (b,u) v 
\\ &= z^{-N} u^{-P} (z-w)^N (u-w)^P
  i_{z,u} i_{w,u} Y^M (Y(a,z-u) b,u-w) v \big|_{w=0} \, ,
\end{split}
\end{equation}
where $N\in\ZZ_+$ is from the associativity axiom and 
$P\in\ZZ_+$ is such that $b_{(k)} v=0$ for $k \geq P$.

Formula \eqref{eq:5.4} implies immediately the following lemma.

\begin{lemma}\label{lem:5.1}
Let $M$ be a module over a field algebra $V$ and let $\A$ be the
associative subalgebra of\/ $\End_{\CC} M$ generated by all
endomorphisms $a^M_{(n)}$ $(a \in V, n \in \ZZ)$.  
Then for any $\al\in\A$, $v\in M$, $\al v$ is a $\CC$-linear
combination of $a^M_{(n)} v$.
\end{lemma}

This lemma along with Jacobson's density theorem
 (see e.g.~\cite{L}) implies:

\begin{theorem}[Density theorem]
  \label{th:5.1}

Let $M$ be a semisimple module over a field algebra $V$ {\rm(}i.e.,~$M$
is a direct sum of irreducible $V$-modules{\rm)}, let $K=\End_V M$, $f
\in \End_K M$ and let $U$ be a finite-dimensional {\rm(}over $\CC)$
subspace of $M$.  Then there exists a finite set of elements $\{
a^i \}_{i \in I}$ of $V$ and a set of integers 
$\{ n_i \}_{i \in I}$ such that
\begin{equation*}
  \sum_i a^{iM}_{(n_i)} |_U = f|_U \, .
\end{equation*}
In particular, if $M$ is an irreducible countable-dimensional
$V$-module, then the conclusion of the theorem holds for any $f
\in \End_{\CC}M$.  

\end{theorem}

Likewise, the proof of Theorem~1.1  from \cite{KR} adapted to the
field algebra framework (by making use of Theorem~\ref{th:5.1} in
place of Jacobson's density theorem) gives the field algebra
analogue of that theorem:

\begin{theorem}[Duality theorem]
  \label{th:5.2}

Let $V$ be a field algebra and let $M$ be an irreducible
$V$-module.  Let $G$ be a group and suppose one has an action of
$G$ on $V$ by automorphisms and a representation of $G$ in $M$
such that the following three properties hold{\rm:}

\begin{list}{}{}
\item \emph{(compatibility)}  
$g (a_{(n)}v) = (g a)_{(n)} (g v)$, 
for $g \in G$, $a \in V$, $v\in M$, $n \in \ZZ;$

\item \emph{(semisimplicity of $G$ on $V$)}  
$V$ is a direct sum of irreducible $G$-modules{\rm;}

\item \emph{(semisimplicity of $G$ on $M$)}
$M$ is a direct sum of at most
countable number of finite-dimensional irreducible $G$-modules.

\end{list}

Then each non-zero isotypic component $M_{\pi}$ {\rm(}where $\pi$ is a
finite-dimensional irreducible $G$-module{\rm)} is irreducible with
respect to the pair $(V^G,G)$ {\rm(}where $V^G$ is the fixed point
subalgebra of the field algebra $V)$.
\end{theorem}

\begin{remark}\label{rem:5.1}
It was stated in \cite{KT} that the proof of Theorem~1.1 from \cite{KR}
can be adapted to give \thref{th:5.2} for vertex algebras. A proof in the
case of graded vertex algebras with the use of Zhu's algebra was given
in \cite{DLM, DY}.
\end{remark}

Consider again formula \eqref{eq:5.4},
where $N\in\ZZ_+$ is from the associativity axiom and 
$P\in\ZZ_+$ is such that $b_{(k)} v=0$ for $k \geq P$.
Taking $\Res_u u^{P-1}$, we obtain:
\begin{equation}\label{eq:5.5}
  Y^M (a,z) (b^M_{(P-1)} v)
  = \sum_{i=0}^N \binom N i z^{-i} (-w)^{P+i} \,
      Y^M (Y(a,z) b, -w) v \big|_{w=0} \, ,
\end{equation}
which implies
\begin{equation}\label{eq:5.6}
a^M_{(n)} (b^M_{(P-1)} v)
  = \sum_{i=0}^N \binom N i (a_{(n-i)} b)^M_{(P+i-1)} v \, ,
\qquad n\in\ZZ \, .
\end{equation}
%
Since the right-hand side of \eqref{eq:5.5}
contains only non-negative powers of $w$, we have:
\begin{equation}\label{eq:5.7}
0 = \sum_{i=0}^N \binom N i (a_{(n-i)} b)^M_{(P+i)} v \, ,
\qquad n\in\ZZ \, .
\end{equation}

\begin{definition}
  \label{def:gm}
A module $M$ over a graded field algebra $(V,\vac,Y)$ is called \emph{graded}
if there is a diagonalizable operator $H^M \in \End M$ satisfying
{\rm(}cf.\ \eqref{eq:g.1}{}{\rm)}
\begin{equation}\label{eq:gm.1}
[H^M, Y^M(a,z)] = Y^M(Ha,z) + z\d_z Y^M(a,z) \, ,
\qquad a \in V \, .
\end{equation}
It is called $\ZZ_+$-\emph{graded} if 
all eigenvalues of $H^M$ are non-negative integers.
\end{definition}

{}From now until the end of the section, 
$(V,\vac,Y)$ will be a $\ZZ_+$-graded field algebra
and $M$ will be a $\ZZ_+$-graded module over it.
Then $M = \bigoplus_{k\in\ZZ_+} M_k$, where
$M_k = \{ v\in M \st H^M v=k v \}$
(so $M_k=0$ for $k\not\in\ZZ_+$).
By \eqref{eq:g.4}, we have:
\begin{equation*}\label{eq:gm.2a}
a^M_n M_k  \subset M_{k-n} \, ,
\qquad a \in V, \; n,k \in \ZZ_+ \, ,
\end{equation*}
where $a^M_n$ is defined similarly to \eqref{eq:g.3}.
In particular,
\begin{equation*}\label{eq:gm.2}
a^M_0 v \in M_0, \quad
a^M_n v = 0 \qquad \text{for all \; $a \in V$, $v \in M_0$, $n>0$.}
\end{equation*}
%

For a homogeneous $a \in V$, $b \in V$, and $r \in \ZZ$, let us define 
\begin{equation}\label{eq:z.1}
a *_r b = \Res_z z^r (z+1)^{\De_a} \, Y(a,z) b
        = \sum_{i=0}^{\De_a} \binom {\De_a} i a_{(r+i)} b \, ,
\end{equation}
and extend it to a bilinear product on $V$.
For $r=-1$, this coincides with Zhu's $*$-product,
and for $r=-2$ with Zhu's $\circ$-product \cite{Z}.
We will also denote $*_{-1}$ as $*$ for short. 
Using $[H,T]=T$, it is easy to see that
\begin{equation}\label{eq:z.2}
((T+H)a) *_r b = -r \, a *_{r-1} b - (r+1) \, a *_r b \, ,
\qquad a,b \in V, \; r \in \ZZ \, .
\end{equation}
In particular, $a *_{-2} b = ((T+H)a) * b$.

\begin{proposition}\label{prop:z.1}
Assume that the elements $a,b \in V$, $v \in M_0$
satisfy the associativity relation from \deref{def:5.1}
with $N=\De_a$. Then{\rm:}
\begin{align}
\label{eq:z.3}
a^M_0 (b^M_0 v) &= (a*b)^M_0 v \, , 
\\
\label{eq:z.4}
0 &= (a *_r b)^M_0 v  \qquad \text{for } \; r \le -2 \, .
\end{align}
\end{proposition}
\begin{proof}
Equation~\eqref{eq:z.3} follows immediately from 
\eqref{eq:g.2a}, \eqref{eq:g.3} and \eqref{eq:5.6}
for $n=\De_a-1$, $N=\De_a$, $P=\De_b$.
Equation~\eqref{eq:z.4} follows from
\eqref{eq:z.2}, \eqref{eq:z.3} and \eqref{eq:g.5}.
\end{proof}


It follows from \reref{rem:g.2} that the assumption of \prref{prop:z.1}
is satisfied when $V$ is a $\ZZ_+$-graded strong field algebra
and $M$ is a $\ZZ_+$-graded module over it.

\begin{proposition}\label{prop:z.2}
If $V$ is a $\ZZ_+$-graded strong field algebra, 
the space $V *_{-2} V = \Span \{ a *_{-2} b \st a,b \in V \}$
is a two-sided ideal with respect to the $*$-product,
and the $*$-product becomes associative on the quotient.
\end{proposition}

The associative algebra $\Zh(V) = V / (V *_{-2} V)$ 
is called the \emph{Zhu algebra} of $V$ (see~\cite{Z}). 
Then, by \prref{prop:z.1}, $M_0$ is a module over $\Zh(V)$.

\begin{proof}[Proof of \prref{prop:z.2}]
First of all, notice that by \eqref{eq:z.2}
\begin{equation}\label{eq:z.5}
a *_r b \in V *_{-2} V \qquad \text{for all $r\le -2$ and $a,b \in V$}.
\end{equation}

We are going to show that for any $a,b,c \in V$
\begin{equation}\label{eq:z.6}
(a *_r b) *_s c \equiv a *_r (b *_s c) \mod V *_{-2} V 
\qquad \text{for $r,s \le -1$}.
\end{equation}
We may assume that $a$ and $b$ are homogeneous.
Then, using \eqref{eq:g.2a},
we find:
\begin{align*}
(a *_r b) *_s c
&= \sum_{i=0}^{\De_a} \binom {\De_a} i 
         \Res_u u^s (u+1)^{\De_a +\De_b-r-i-1} \, Y(a_{(r+i)} b, u)c
\\
&= \Res_z \Res_u z^r u^s (u+1)^{\De_b-r-1} (z+u+1)^{\De_a} \, 
   Y(Y(a,z)b,u)c \, .
\end{align*}
Now we use \eqref{eq:3.8} with $w=-u$.
Making the change of variables $x=z+u$, we obtain
for any $b' \in V$:
\begin{align*}
\Res_z & z^r (z+u+1)^{\De_a} \, i_{z,u} \, Y(a,z+u) b'
= \Res_x  i_{x,u} (x-u)^r (x+1)^{\De_a} \, Y(a,x) b'
\\
&= \sum_{k=0}^\infty \binom r k (-u)^k \, a *_{r-k} b'
\equiv a *_r b' \mod V *_{-2} V \, 
\qquad \text{for $r\le -1$}.
\end{align*}
%
Therefore, the first term in the right-hand side of \eqref{eq:3.8}
gives $a *_r (b *_s c)$ modulo $V *_{-2} V$.
The contribution of the second term is zero, because for 
$r\le -1$ and $j\ge0$
\begin{equation*}
\Res_z z^r (z+u+1)^{\De_a} \, \d_u^j \de(z-(-u))
\in u^r \CC[u^{-1}] {} \, .
\end{equation*}
This completes the proof of \eqref{eq:z.6}, and hence of
the proposition.
\end{proof}

A large part of Zhu's theory can be generalized from vertex algebras
to the case of strong field algebras.
For example, we have the following result (cf.\ \cite{Z}).

\begin{proposition}\label{prop:z.3}
If $V$ is a $\ZZ_+$-graded strong field algebra,
the functor $M \mapsto M_0$ establishes a one-to-one
correspondence between irreducible $\ZZ_+$-graded $V$-modules 
and irreducible $\Zh(V)$-modules.
\end{proposition}

However, because of \thref{thm:sf.1}, the case of strong field algebras
is not much different than the case of vertex algebras.
Some of the results hold in the more interesting
case of field algebras (see e.g.~\eqref{eq:5.6}, \eqref{eq:5.7}),
but we do not know whether $V *_{-2} V$ is an ideal or
$ V / (V *_{-2} V)$ is associative in such generality.


\begin{example}\label{ex:z4}
Let $V$ be the tensor product of a $\ZZ_+$-graded strong field algebra
$V'$ and a unital algebra $A$ (see \exref{ex:1.4}).
Then, although $V$ is not a strong field algebra, 
Propositions~\ref{prop:z.1}--\ref{prop:z.3} hold for $V$. 
One has: $\Zh(V) \simeq \Zh(V') \tt A$.
\end{example}

\begin{example}\label{ex:z5}
Similarly, let $V = V' \smash \Ga$ be the smash product of 
a $\ZZ_+$-graded strong field algebra $V'$ and a group $\Ga$ of 
its automorphisms (see \exref{ex:1.5}). 
Then Propositions~\ref{prop:z.1}--\ref{prop:z.3} hold for $V$. 
Note that the action of $\Ga$ on
$V'$ induces an action on its Zhu algebra $\Zh(V')$.
One has: $\Zh(V) \simeq \Zh(V') \smash \Ga$.
\end{example}

%



\section*{Acknowledgments}

We are grateful to MSRI, where this work was completed, for hospitality.
One of the authors wishes to thank V.~Schechtman and L.~W.~Small
for very useful discussions and V.~L.~Popov and E.~B.~Vinberg 
for correspondence.
We thank A.~De~Sole for carefully reading
the manuscript and suggesting improvements,
and C.~Dong for pointing out a gap in a previous version
of the paper.


\end{document}